\documentclass[12pt]{article}

\usepackage{amsmath}
\usepackage{amsfonts}
\usepackage{amssymb}
\usepackage{theorem}

\theorembodyfont{\it}
\newtheorem{theorem}{Theorem}[section]
\newtheorem{corollary}[theorem]{Corollary}

\newtheorem{lemma}[theorem]{Lemma}

\theorembodyfont{\rm}

\newtheorem{notation}[theorem]{Notation}

\def\proof{\noindent{\bf Proof: }}
\def\endproof{\hfill$\square$}

\renewenvironment{enumerate}
{
    \begin{list}{\makebox[1.0cm]{\emph{(\roman{enumi})}}}
    {
        \usecounter{enumi}
        \setlength{\topsep}{0pt}
        \setlength{\partopsep}{0pt}
        \setlength{\parsep}{0pt}
        \setlength{\itemsep}{0pt}
        \setlength{\labelsep}{0in}
        \setlength{\labelwidth}{1.0cm}
        \setlength{\leftmargin}{1.0cm}
    }
} {
    \end{list}
}

\begin{document}

\thispagestyle{empty}

\title{Contact Systems and Corank One \\ Involutive Subdistributions}

\author{William {\sc Pasillas-L\'epine}\thanks{Now at Laboratoire des signaux et syst\`emes,
CNRS -- Sup\'elec. Plateau de Moulon, 3 rue Joliot-Curie. 91 192 Gif-sur-Yvette
Cedex, France. E-mail: \texttt{pasillas@lss.supelec.fr}} \hbox{ and }Witold
{\sc Respondek} \bigskip \\ Institut national des sciences appliqu\'ees de
Rouen\\ D\'epartement g\'enie math\'ematique\\ Place \'Emile Blondel \\ 76 131
Mont Saint Aignan Cedex, France \bigskip \\ Tel: (+33) 02 35 52 84 32 ; Fax:
(+33) 02 35 52 83 32 \\ E-Mail: {\tt wresp@lmi.insa-rouen.fr}}

\date{March, 2000}

\enlargethispage{2cm}

\maketitle

\begin{abstract}
We give necessary and sufficient geometric conditions for a distribution (or a
Pfaffian system) to be locally equivalent to the canonical contact system on
$J^{n}(\mathbb{R},\mathbb{R}^{m})$. We study the geometry of that class of
systems, in particular, the existence of corank one involutive
subdistributions. We also distinguish regular points, at which the system is
equivalent to the canonical contact system, and singular points, at which we
propose a new normal form that generalizes the canonical contact system on
$J^{n}(\mathbb{R},\mathbb{R}^{m})$ in a way analogous to that how Kumpera-Ruiz
normal form generalizes the canonical contact system on
$J^{n}(\mathbb{R},\mathbb{R})$, which is also called Goursat normal form.

\medskip


\medskip


\end{abstract}

\newpage

\section*{Introduction}

Consider $J^{n}(\mathbb{R}^{k},\mathbb{R}^{m})$, the space of $n$-jets of
smooth maps from $\mathbb{R}^{k}$ into $\mathbb{R}^{m}$, and denote by
\[
(q^{1},\ldots,q^{k},u_{1},\ldots,u_{m},p_{i}^{\sigma}), \text{\quad for }1\leq
i\leq m\text{ and for }1\leq \left| \sigma \right| \leq n \text{,}
\]
the \emph{canonical coordinates}, also called natural coordinates, on this
space (see e.g. \cite{bryant-chern-gardner-goldschmidt-griffiths},
\cite{olver-equivalence}, and \cite{vinogradov}), where $q^{j}$, for $1\leq
j\leq k$, represent \emph{independent variables} and $u_{i}$, for $1\leq i\leq
m$, represent \emph{dependent variables}; the vector of non-negative integers
$\sigma=(\sigma_{1},\ldots,\sigma_{k}) $ is a multi-index such that $\left|
\sigma\right|  =\sigma_{1}+\cdots+\sigma _{k}\leq n$; and $p_{i}^{\sigma}$, for
$1\leq i\leq m$,
correspond to the partial derivatives $\partial^{\left|  \sigma\right|  }%
u_{i}/\partial q_{\sigma}$. Any smooth map~$\varphi$ from~$\mathbb{R}^{k}$
into~$\mathbb{R}^{m}$ defines a submanifold in $J^{n}(\mathbb{R}%
^{k},\mathbb{R}^{m})$ by the relations
\[
p_{i}^{\sigma}=\frac{\partial^{\left|  \sigma\right|  }\varphi_{i}}{\partial
q_{\sigma}}(q^{1},\ldots,q^{k}),
\]
for $1\leq i\leq m$ and $1\leq\left|  \sigma\right|  \leq n$. This submanifold
is called the $n$\emph{-graph} of~$\varphi$. It turns out that all $n$-graphs
are integral submanifolds, of dimension~$k$, of a distribution called the
\emph{canonical contact system} on $J^{n}(\mathbb{R}^{k},\mathbb{R}^{m})$ or
the \emph{Cartan distribution} \cite{vinogradov} on
$J^{n}(\mathbb{R}^{k},\mathbb{R}^{m})$. The Pfaffian system that anihilates
this distribution, which is also called the \emph{canonical contact system}
(see e.g
\cite{bryant-chern-gardner-goldschmidt-griffiths} and \cite{olver-equivalence}%
), is given in the canonical coordinates of $J^{n}(\mathbb{R}^{k}%
,\mathbb{R}^{m})$ by
\[
dp_{i}^{\sigma}-%
{\textstyle\sum\limits_{j=1}^{k}}
p_{i}^{\sigma+1_{j}}dq^{j}=0,\text{\quad for }1\leq i\leq m\text{ and for
}1\leq \left| \sigma \right| \leq n-1\text{,}
\]
where $\sigma+1_{j}=(\sigma_{1},\ldots,\sigma_{j}+1,\ldots,\sigma_{k})$.


The above description explains the importance of contact systems in geometric
theory of (partial) differential equations and in differential geometry. In the
former, a (partial) differential equation is interpreted as a submanifold in
$J^{n}(\mathbb{R}^{k},\mathbb{R}^{m})$ and thus it is natural to study the
geometry of pairs consisting of a contact system and a submanifold, see
e.g.~\cite{bryant-chern-gardner-goldschmidt-griffiths},
\cite{lychagin-classification}, and~\cite{vinogradov}. In the latter, contact
systems describe, for instance, diffeomorphisms which preserve the $n$-graphs
of appplications (for example, $n$-graphs of curves in the case $k=1$), see
e.g.~\cite{bryant-chern-gardner-goldschmidt-griffiths}
and~\cite{olver-equivalence}.

A natural problem which arises is to characterize those distributions which are
(locally) equivalent to a canonical contact system. This problem was posed by
Pfaff~ \cite{pfaff} in 1815 and seems still to be open in its full generality
although many important particular solutions have been obtained. In the case
$n=1$, with $m=1$ and an arbitrary $k$ the final solution has been obtained by
Darboux~\cite{darboux} in his famous theorem generalizing earlier results of
Pfaff~\cite{pfaff} and Frobenius~\cite{frobenius}. The case $n=2$, $m=1$ and
$k=1$ was solved by Engel in~\cite{engel}. The case $n\geq2$, $m=1$ and $k=1$
was solved by E. von Weber~\cite{weber-article},
Cartan~\cite{cartan-equivalence-absolue} and Goursat~\cite{goursat} (at generic
points) and by Libermann~\cite{libermann}, Kumpera and
Ruiz~\cite{kumpera-ruiz}, and Murray~\cite{murray-nilpotent} (at an arbitrary
point). The case $n=1$, with $k$ and $m$ arbitrary has been studied and solved
by Bryant \cite{bryant-thesis} (see also
\cite{bryant-chern-gardner-goldschmidt-griffiths}).

\medskip

This paper is devoted to the problem of when a given distribution is locally
equivalent to the canonical contact system in the case $k=1$, $n$ and $m$
arbitrary, that is to the canonical contact system for curves. This problem has
been studied by Gardner and Shadwick \cite{gardner-shadwick},
Murray~\cite{murray-nilpotent}, and Tilbury and Sastry~\cite{tilbury-sastry}
(as a particular case of the the problem of equivalence to the so-called
extended Goursat normal form). Their solutions are based on a result of
\cite{gardner-shadwick} that assures the equivalence provided that a certain
differential form satisfies precise congruence relations. The problem of how to
verify the existence of such a form had apparently remained open. This
difficulty was solved by Aranda-Bricaire and Pomet~\cite{aranda-pomet-nolcos},
who proposed an algorithm which determines the existence of such a form. Their
solution, although being elegant and checkable, uses the formalism of infinite
dimensional manifolds and thus goes away from classical results characterizing
contact systems.

The aim of our paper is two-fold. Firstly, in Theorem ~\ref{thm-pfaff-curves},
we give geometric checkable conditions, based on the classical notion of
Engel's rank, which characterize regular contact systems for curves, i.e., for
$k=1$ and arbitrary $n$ and $m$. Secondly, we extend our approach to singular
points and we prove, in Theorem~\ref{thm-ext-kumpera-ruiz}, that any singular
contact system can be put to a normal form for which we propose the name
extended Kumpera-Ruiz normal form. That form generalizes canonical contact
systems for curves in a way analogous to that how Kumpera-Ruiz normal form
(see~\cite{kumpera-ruiz}; compare \cite{cheaito-mormul},
\cite{cheaito-mormul-pasillas-respondek}, \cite{montgomery-zhitomirskii},
\cite{mormul-codimension-one}, and~\cite{pasillas-respondek-goursat})
generalizes Goursat normal form, that is the canonical contact system on
$J^{n}(\mathbb{R},\mathbb{R})$. When checking conditions of
Theorem~\ref{thm-pfaff-curves}
 we must determine whether a given distribution possesses a corank one
involutive subdistribution. An elegant answer given to this problem by
Bryant~\cite{bryant-thesis} implies that our conditions become checkable.

\medskip

The paper is organized as follows. In
Section~\ref{sec-canonical-contact-system} we define the canonical system for
curves and we give the first main result of the paper, Theorem
\ref{thm-pfaff-curves}, characterizing distributions that are locally
equivalent to a canonical contact system for curves. In
Section~\ref{sec-corank-one} we discuss the problem of whether a given
distribution possesses a corank one involutive subdistribution and we recall
Bryant's solution of this problem (see also Appendices~\ref{app-A}
and~\ref{app-B}). In Section~\ref{sec-extended-kumpera-ruiz} we introduce
extended Kumpera-Ruiz normal forms and give the second main result of the
paper, Theorem~\ref{thm-ext-kumpera-ruiz}, stating that any singular contact
system is locally equivalent to an extended Kumpera-Ruiz normal form. Finally,
Section~\ref{sec-two-main} contains proofs of all our results.

\medskip

\noindent {\bf Acknowledgments:} The authors are grateful to Robert Bryant, who
kindly sent them a copy of his Ph.D. Thesis~\cite{bryant-thesis} and to
Jean-Baptiste Pomet for discussion on his results~\cite{aranda-pomet-nolcos}.

\section{The Canonical Contact System for Curves}

\label{sec-canonical-contact-system}

A rank\emph{\ }$k$ \emph{distribution }$\mathcal{D}$ on a smooth manifold~$M$
is a map that assigns smoothly to each point~$p$ in$~M$ a linear subspace
$\mathcal{D}(p)\subset T_{p}M$ of dimension~$k$. In other words, a rank~$k$
distribution is a smooth rank~$k$ subbundle of the tangent bundle~$TM$. Such a
field of tangent~$k$-planes is spanned locally by~$k$ pointwise linearly
independent smooth vector fields $f_{1},\ldots,f_{k}$ on~$M$, which will be
denoted by $\mathcal{D}=(f_{1},\ldots,f_{k})$. Two distributions~$\mathcal{D}$
and~$\tilde{\mathcal{D}}$ defined on two manifolds~$M$ and~$\tilde M$,
respectively, are \emph{equivalent} if there exists a smooth
diffeomorphism~$\varphi$ between~$M$ and~$\tilde M$ such that $(\varphi
_{*}\mathcal{D})(\tilde p)=\tilde{\mathcal{D}}(\tilde p)$, for each
point$~\tilde p$ in~$\tilde M$.

The \emph{derived flag} of a distribution~$\mathcal{D}$ is the sequence of
modules of vector fields $\mathcal{D}^{(0)}\subset\mathcal{D}^{(1)}%
\subset\cdots$ defined inductively by
\begin{equation}
\mathcal{D}^{(0)}=\mathcal{D}\text{\quad and\quad}\mathcal{D}^{(i+1)}%
=\mathcal{D}^{(i)}+[\mathcal{D}^{(i)},\mathcal{D}^{(i)}]\text{, \quad for
}i\geq0\text{.}\label{acta-derived-flag}%
\end{equation}
The \emph{Lie flag} is the sequence of modules of vector fields $\mathcal{D}%
_{0}\subset\mathcal{D}_{1}\subset\cdots$ defined inductively by
\begin{equation}
\mathcal{D}_{0}=\mathcal{D}\text{\quad and\quad}\mathcal{D}_{i+1}%
=\mathcal{D}_{i}+[\mathcal{D}_{0},\mathcal{D}_{i}]\text{, \quad for }%
i\geq0\text{.}\label{acta-lie-flag}%
\end{equation}
In general, the derived and Lie flags are different; though for any point~$p$
in the underlying manifold the inclusion $\mathcal{D}_{i}(p)\subset
\mathcal{D}^{(i)}(p)$ clearly holds, for~$i\geq0$.

For a given distribution $\mathcal{D}$, defined on a manifold~$M$, we will say
that a point~$p$ of~$M$ is a \emph{regular point} of~$\mathcal{D}$ if all the
elements $\mathcal{D}_{i}$ of its Lie flag have constant rank in a small enough
neighborhood of $p$. A distribution~$\mathcal{D}$ is said to be
\emph{completely nonholonomic} if, for each point~$p$ in~$M$, there exists an
integer~$N(p)$ such that $\mathcal{D}_{N(p)}(p)=T_{p}M$. A
distribution~$\mathcal{D}$ is said to be \emph{involutive} if its first derived
system satisfies~$\mathcal{D}^{(1)}=\mathcal{D}^{(0)}$.

\medskip

An alternative description of the above defined objects can also be given using
the dual language of differential forms. A \emph{Pfaffian system} $\mathcal{I}$
of rank~$s$ on a smooth manifold~$M$ is a map that assigns smoothly to each
point~$p$ in$~M$ a linear subspace $\mathcal{I}(p)\subset T_{p}^{\ast }M$ of
dimension~$s$. In other words, a Pfaffian system of rank~$s$ is a smooth
subbundle of rank~$s$ of the cotangent bundle~$T^{\ast}M$. Such a field of
cotangent~$s$-planes is spanned locally by~$s$ pointwise linearly independent
smooth differential $1$-forms $\omega_{1},\ldots,\omega_{s}$ on~$M$, which will
be denoted by $\mathcal{I}=(\omega_{1},\ldots,\omega_{s})$. Two Pfaffian
systems~$\mathcal{I}$ and~$\tilde{\mathcal{I}}$ defined on two manifolds~$M$
and~$\tilde{M}$, respectively, are \emph{equivalent} if there exists a smooth
diffeomorphism~$\varphi$ between~$M$ and~$\tilde{M}$ such that $\mathcal{I}%
(p)=(\varphi^{\ast}\tilde{\mathcal{I}})(p)$, for each point$~p$ in~$M$.

For a Pfaffian system $\mathcal{I}$, we can define its \emph{derived flag}
$\mathcal{I}^{(0)}\supset\mathcal{I}^{(1)}\supset\cdots$ by the relations
$\mathcal{I}^{(0)}=\mathcal{I}$ and $\mathcal{I}^{(i+1)}=\{\alpha
\in\mathcal{I}^{(i)}:d\alpha\equiv0\mod\mathcal{I}^{(i)}\}$, for $i \geq 0$,
provided that each element $\mathcal{I}^{(i)}$ of this sequence has constant
rank. In this case, it is immediate to see that the derived flag of the
distribution $\mathcal{D}=\mathcal{I}^{\perp}$ coincides with the sequence of
distributions that anihilate the elements of the derived flag of $\mathcal{I}$,
that is $$ \mathcal{D}^{(i)}=(\mathcal{I}^{(i)})^\perp,\mbox{\quad for } i \geq
0. $$ For a given Pfaffian system $\mathcal{I}$, we will say that a point $p$
of $M$ is a \emph{regular point} if~$p$ is a regular point for the distribution
$\mathcal{D}=\mathcal{I}^{\perp}$, that is if all elements $\mathcal{D}_i$ of
the Lie flag are of constant rank in a small enough neighborhood of~$p$.

\medskip

Consider the space $J^{n}(\mathbb{R},\mathbb{R}^{m})$ of jets of order~$n\geq1$
of functions from~$\mathbb{R}$ into~$\mathbb{R}^{m}$. This space is
diffeomorphic to~$\mathbb{R}^{(n+1)m+1}$. The canonical coordinates associated
to~$\mathbb{R}$ (denoted by $x_{0}^{0}$) and to $\mathbb{R}^{m}$ (denoted by
$x_{1}^{0},\ldots,x_{m}^{0}$) can be used to define the \emph{canonical
coordinates} on $J^{n}(\mathbb{R},\mathbb{R}^{m})$, which will be denoted by
\[
x_{0}^{0},x_{1}^{0},\ldots,x_{m}^{0},x_{1}^{1},\ldots,x_{m}^{1},\ldots
,x_{1}^{n},\ldots,x_{m}^{n},
\]
with obvious indentifications $x_{0}^{0}=q$ and $x_{i}^{0}=u_{i}$, for $1\leq
i\leq m$, and $x_{i}^{j}=p_{i}^{j}$, for $1\leq i\leq m$ and  $1\leq j \leq n$
(see the beginning of the Introduction). Observe that any smooth map~$\varphi$
from~$\mathbb{R}$ into~$\mathbb{R}^{m}$
defines a curve in $J^{n}(\mathbb{R},\mathbb{R}^{m})$ by the relations~$x_{j}%
^{i}=\varphi_{j}^{(i)}(x_{0}^{0})$, for $0\leq i\leq n$ and $1\leq j\leq m$,
where $\varphi_{j}^{(i)}$ denotes the $i$-th derivative with respect to
$x_{0}^{0}$ of the $j$-th component of $\varphi$. This curve is called the
$n$-\emph{graph} of~$\varphi$. It is clear that not all curves in
$J^{n}(\mathbb{R},\mathbb{R}^{m})$ are $n$-graphs of maps. In order to
distinguish the ``good'' curves from the ``bad'' ones, we should introduce a
set of constraints on the velocities of curves in
$J^{n}(\mathbb{R},\mathbb{R}^{m})$. In other words, we should endow
$J^{n}(\mathbb{R},\mathbb{R}^{m})$ with a nonholonomic structure.

The \emph{canonical contact system} on $J^{n}(\mathbb{R},\mathbb{R}^{m})$ is
the completely nonholonomic distribution spanned by the following family of
vector fields:
\begin{equation}
\left(
\begin{array}
[c]{l}%
\tfrac\partial{\partial x_{1}^{n}}%
\end{array}
,\ldots,
\begin{array}
[c]{l}%
\tfrac\partial{\partial x_{m}^{n}}%
\end{array}
,
\begin{array}
[c]{l}%
{\textstyle\sum\limits_{i=0}^{n-1}}
{\textstyle\sum\limits_{j=1}^{m}}
x_{j}^{i+1}\tfrac\partial{\partial x_{j}^{i}}+\tfrac\partial{\partial
x_{0}^{0}}%
\end{array}
\right)  .\label{contact-system}%
\end{equation}
By definition, if a curve in $J^{n}(\mathbb{R},\mathbb{R}^{m})$ is the
$n$-graph of some map then it is an integral curve of the canonical contact
system. More precisely, a section $\sigma : \mathbb{R} \rightarrow
J^{n}(\mathbb{R},\mathbb{R}^{m})$ is the $n$-graph of a curve $\varphi :
\mathbb{R} \rightarrow \mathbb{R}^n$ if and only if it is an integral curve of
the canonical contact system on $J^{n}(\mathbb{R},\mathbb{R}^{m})$ (see e.g.
\cite{olver-equivalence}).

\medskip

\noindent The aim of our paper is to give a complete answer to the question
``Which distributions are locally equivalent to the canonical contact system on
$J^{n}(\mathbb{R},\mathbb{R}^{m})$?''.\ The following result will be the
starting point of our study.

\begin{theorem}[contact systems for curves]
\label{thm-pfaff-curves}A rank $m+1$ distribution~$\mathcal{D}$ on a
manifold~$M$ of dimension~$(n+1)m+1$ is equivalent, in a small enough
neighborhood of any point~$p$ in~$M$, to the canonical contact system on
$J^n(\mathbb{R},\mathbb{R}^m)$ if and only if the two following conditions
hold, for $0\leq i\leq n$.
\begin{enumerate}
\item Each element~$\mathcal{D}^{(i)}$ of the derived flag has constant
rank~$(i+1)m+1$ and contains an involutive subdistribution~$\mathcal{L}%
_i\subset\mathcal{D}^{(i)}$ that has constant corank one in~$\mathcal{D}%
^{(i)}$.
\item Each element~$\mathcal{D}_i$ of the Lie flag has constant rank~$%
(i+1)m+1$.
\end{enumerate}
\end{theorem}

\noindent This result yields a constructive test for the local equivalence to
the canonical contact system for curves, provided that we know how to check
whether or not a given distribution admits a corank one involutive
subdistribution. We give in the next section a checkable necessary and
sufficient condition for the existence of such a distribution. The proof of
Theorem~\ref{thm-pfaff-curves} will be given in Section~\ref{sec-two-main}.

\section{Corank One Involutive Subdistributions}

\label{sec-corank-one}

The aim of this Section is to give an answer to the following question: ``When
does a given constant rank distribution $\mathcal{D}$ contain an involutive
subdistribution $\mathcal{L}\subset\mathcal{D}$ that has constant corank one in
$\mathcal{D}$?''. In fact, the answer to this question is an immediate
consequence of a result contained in Bryant's Ph.D.
thesis~\cite{bryant-thesis}. Links between Bryant's result and the
characterization of the canonical contact system for curves have also been
observed by Aranda-Bricaire and Pomet~\cite{aranda-pomet-nolcos}.

A \emph{characteristic vector field} of a distribution $\mathcal{D}$ is a
vector field $f$ that belongs to $\mathcal{D}$ and satisfies $[f,\mathcal{D}%
]\subset\mathcal{D}$. The \emph{characteristic distribution} of $\mathcal{D}$,
which will be denoted by $\mathcal{C}$, is the module spanned by all its
characteristic vector fields. It follows directly from the Jacobi identity that
the characteristic distribution is always involutive. For a constant rank
Pfaffian system $\mathcal{I}$, the characteristic distribution of
$\mathcal{I}^{\perp}$ is often called the \emph{Cartan system} of
$\mathcal{I}$. We refer the reader
to~\cite{bryant-chern-gardner-goldschmidt-griffiths} for a definition of the
Cartan system given in the language of Pfaffian systems.

The \emph{Engel rank} \cite{bryant-chern-gardner-goldschmidt-griffiths} of a
Pfaffian system $\mathcal{I}$, at a point~$p$, is the largest integer~$\rho$
such that there exists a $1$-form $\omega$ in $\mathcal{I}$ for which we have
$(d\omega)^{\rho}(p)\neq0\operatorname*{mod}\mathcal{I}$. The Engel rank of a
constant rank distribution $\mathcal{D}$ will be, by definition, the Engel rank
of its anihilator$~\mathcal{D}^{\perp}$. Obviously, the Engel rank~$\rho$ of a
distribution equals zero at each point if and only if the distribution is
involutive.

We will now give an equivalent definition of the Engel rank in the language of
vector fields, in the particular case when $\rho=1$, which will be important in
the paper. Let $\mathcal{D}$ be a distribution such that $\mathcal{D}^{(0)}$
and $\mathcal{D}^{(1)}$ have constant ranks $d_{0}$ and $d_{1}$, respectively,
and denote $r_{0}(p)=d_{1}(p)-d_{0}(p)$. Assume that $d_{0}\geq2$ and
$r_{0}\geq1$. Take a family of vector fields
\[
(f_{1},\ldots,f_{d_{0}},g_{1},\ldots,g_{r_{0}})
\]
such that $\mathcal{D}^{(0)}=(f_{1},\ldots,f_{d_{0}})$ and$~\mathcal{D}%
^{(1)}=(f_{1},\ldots,f_{d_{0}},g_{1},\ldots,g_{r_{0}})$. The \emph{structure
functions} $c_{ij}^{k}$ associated to those generators, for $1\leq i<j\leq
d_{0}$ and $1\leq k\leq r_{0}$, are the smooth functions defined by the
following relations:
\[
\lbrack f_{i},f_{j}]=%
{\textstyle\sum\limits_{k=1}^{r_{0}}}
c_{ij}^{k}\,g_{k}\operatorname*{mod}\mathcal{D}^{(0)},\text{\quad for }1\leq
i<j\leq d_{0}.
\]
It is important to point out that the structure functions are not invariantly
related to the distribution $\mathcal{D}$, since they depend on the choice of
generators.

Assume that the Engel rank $\rho$ of $\mathcal{D}$ is constant and that
$r_{0}\geq1$. It is easy to check that $\rho=1$ if and only if either
$d_{0}=2$, or $d_{0}=3$, or $d_{0}\geq4$ and the structure functions satisfy
the relations
\begin{equation}
c_{ij}^{p}c_{kl}^{q}-c_{ik}^{p}c_{jl}^{q}+c_{il}^{p}c_{jk}^{q}+c_{jk}%
^{p}c_{il}^{q}-c_{jl}^{p}c_{ik}^{q}+c_{kl}^{p}c_{ij}^{q}=0,\label{engel-rank}%
\end{equation}
for each sextuple $(i,j,k,l,p,q)$ of integers such that $1\leq i<j<k<l\leq
d_{0}$ and $1\leq p\leq r_{0}$ and $1\leq q\leq r_{0}$.

The following result is a particular case of Bryant's algebraic
Lemma~\cite{bryant-thesis} (see
also~\cite{bryant-chern-gardner-goldschmidt-griffiths}), which is the
cornerstone of Bryant's characterization of the canonical contact system on
$J^{1}(\mathbb{R}^{k},\mathbb{R}^{m})$.

\begin{lemma}[Bryant]
Let $\mathcal{D}$ be a distribution such that $\mathcal{D}^{(0)}$ and $%
\mathcal{D}^{(1)}$ have constant ranks $d_0$ and $d_1$, respectively. Assume
that $r_0\geq1$.\ Then the two following conditions are equivalent:
\begin{enumerate}
\item The characteristic distribution $\mathcal{C}$ of $\mathcal{D}$ has
constant rank $c_0=d_0-r_0-1$ and the Engel rank $\rho$ of $\mathcal{D}$ is
constant and equals$~1$;
\item The distribution $\mathcal{D}$ contains a subdistribution $\mathcal{B}%
\subset\mathcal{D}$ that has constant corank one in $\mathcal{D}$ and satisfies
$[\mathcal{B},\mathcal{B}]\subset\mathcal{D}$.
\end{enumerate}
\end{lemma}

Observe that if the first condition is satisfied then we must necessarily have
$r_{0}\leq d_{0}-1$. The following result is included in the proof of
Bryant's~\cite{bryant-thesis} normal form Theorem (see
also~\cite{bryant-chern-gardner-goldschmidt-griffiths}). In
Appendix~\ref{app-A} we give an alternative proof of its surprising Item (iii)
; our proof explains the role of the assumption $r_0 \geq 3$ by relating it to
the Jacobi identity.

\begin{lemma}[Bryant]
Let $\mathcal{D}$ be a distribution such that $\mathcal{D}^{(0)}$ and $%
\mathcal{D}^{(1)}$ have constant ranks $d_0$ and $d_1$, respectively. Assume
that the distribution $\mathcal{D}$ contains a subdistribution $\mathcal{B}%
\subset\mathcal{D}$ that has constant corank one in $\mathcal{D}$ and satisfies
$[\mathcal{B},\mathcal{B}]\subset\mathcal{D}$.
\begin{enumerate}
\item If $r_0=1$ then the distribution $\mathcal{D}$ contains an involutive subdistribution $%
\mathcal{L}\subset\mathcal{D}$ that has constant corank one in $\mathcal{D}$;
\item If $r_0\geq2$ then $\mathcal{B}$ is unique;
\item If $r_0\geq3$ then $\mathcal{B}$ is involutive.
\end{enumerate}
\end{lemma}

Observe that, in the first item of the above Lemma, the involutive
subdistribution $\mathcal{L}$ can be different from $\mathcal{B}$, which is not
necessarily involutive. The following result is a direct consequence of
Bryant's work.\ In order to avoid the trivial case $r_{0}=0$, for which the
existence of a corank one involutive subdistribution is obvious, we will assume
that $r_{0}\geq1$.

\begin{corollary}[corank one involutive subdistributions]
\label{cor-corank-one-involutive} Let $\mathcal{D}$ be a distribution such that
$\mathcal{D}^{(0)}$ and $\mathcal{D}^{(1)}$ have constant ranks $d_0$
and $d_1$, respectively. Assume that $r_0\geq1$. Then, the distribution $%
\mathcal{D}$ contains an involutive subdistribution $\mathcal{L}\subset
\mathcal{D}$ that has constant corank one in$~\mathcal{D}$ if and only if the
three following conditions hold:
\begin{enumerate}
\item The characteristic distribution $\mathcal{C}$ of $\mathcal{D}$ has
constant rank $c_0=d_0-r_0-1$;
\item The Engel rank $\rho$ of $\mathcal{D}$ is constant and equals$~1$;
\item If $r_0=2$ then, additionally, the unique corank one subdistribution $%
\mathcal{B}\subset\mathcal{D}$ such that $[\mathcal{B},\mathcal{B}]\subset
\mathcal{D}$ must be involutive.
\end{enumerate}
\end{corollary}


We would like to to emphasize that the above conditions are easy to verify, as
well as the conditions of Corollary~\ref{cor-slide} below. Indeed, for any
distribution, or the corresponding Pfaffian system, we can compute the
characteristic distribution ${\cal C}$ and check whether or not the Engel rank
equals~$1$ using, respectively, the formula~(\ref{cond-C}) and the
condition~(\ref{cond-E}) of Appendix~\ref{app-B}. This gives the solution if
$r_{0}\neq 2$. If $r_{0}=2$ we have additionally to check the involutivness of
the unique distribution ${\cal B}$ satisfying $[{\cal B},{\cal B}]\subset {\cal
D}$, whose explicit construction is also given in Appendix~\ref{app-B}.

Combining Theorem~\ref{thm-pfaff-curves} and
Corollary~\ref{cor-corank-one-involutive} we get the following characterization
of the canonical system on~$J^{n}(\mathbb{R},\mathbb{R}^{m})$, using the
language of Pfaffian systems.

\begin{corollary}[contact systems for curves]
\label{cor-slide} Let $\mathcal{I}$ be a Pfaffian system of rank $nm$, defined
on a manifold~$M$of dimension $(n+1)m+1$. If $m \neq 2$, the Pfaffian
system~$\mathcal{I}$ is locally equivalent, at a given point~$p$ of~$M$, to the
canonical contact system on $J^{n}(\mathbb{R},\mathbb{R}^{m})$ if and only if
\begin{enumerate}
\item The rank of each derived system $\mathcal{I}^{(i)}$ is constant and
      equals $(n-i)m$, for $0\leq i\leq n$;
\item The Engel rank of~$\mathcal{I}^{(i)}$ is constant and equals~$1$,
      for $0\leq i\leq n$;
\item The rank of each Cartan system $\mathcal{C}(\mathcal{I}^{(i)})$ is
      constant and equals $(n+1-i)m+1$, for $0\leq i\leq n-1$ ;
\item The point~$p$ is a regular point for $\mathcal{I}$.
\end{enumerate}
\end{corollary}

\medskip

In other words, if $m \neq 2$, the characterization of the canonical contact
system on~$J^{n}(\mathbb{R},\mathbb{R}^{m})$ turns out to be a natural
combination of that given for $J^{1}(\mathbb{R},\mathbb{R}^{m})$ by Bryant
(see~\cite{bryant-thesis}
and~\cite{bryant-chern-gardner-goldschmidt-griffiths}) and that given for
$J^{n}(\mathbb{R},\mathbb{R})$ by Murray~\cite{murray-nilpotent}.

\section{Extended Kumpera-Ruiz Normal Forms}

\label{sec-extended-kumpera-ruiz}

The aim of this Section is to study the class of distributions that satisfy
condition~(i) of Theorem~\ref{thm-pfaff-curves} but fail to satisfy
condition~(ii) of that theorem. The fist condition describes the geometry of
the canonical contact system while the second condition characterizes regular
points. In this sense, systems that satisfy the former but fail to satisfy the
latter can be considered as ``singular'' contact systems for curves. We will
show that any such distribution can be brought to a normal form for which we
propose the name extended Kumpera-Ruiz normal form.\ Those forms generalize the
canonical contact system on $J^{1}(\mathbb{R},\mathbb{R}^{m})$ in a way
analogous to that how Kumpera-Ruiz normal forms (see e.g. \cite{kumpera-ruiz},
\cite{montgomery-zhitomirskii}, \cite{mormul-codimension-one},
and~\cite{pasillas-respondek-goursat}) generalize the canonical system on
$J^{1}(\mathbb{R},\mathbb{R})$, which is also called Goursat normal form.

Consider the family of vector fields $\kappa^{1}=(\kappa_{1}^{1},\ldots
,\kappa_{m}^{1},\kappa_{0}^{1})$ that span the canonical contact system on
$J^{1}(\mathbb{R},\mathbb{R}^{m})$, where
\begin{align*}
\kappa_{1}^{1}  & =\tfrac\partial{\partial x_{1}^{1}},\ldots,\kappa_{m}%
^{1}=\tfrac\partial{\partial x_{m}^{1}}\\ \kappa_{0}^{1}  &
=x_{1}^{1}\tfrac\partial{\partial x_{1}^{0}}+\cdots
+x_{m}^{1}\tfrac\partial{\partial x_{m}^{0}}+\tfrac\partial{\partial x_{0}%
^{0}},
\end{align*}
and the the family of vector fields $\kappa^{2}=(\kappa_{1}^{2},\ldots
,\kappa_{m}^{2},\kappa_{0}^{2})$ that spans the canonical contact system on
$J^{2}(\mathbb{R},\mathbb{R}^{m})$, where
\begin{align*}
\kappa_{1}^{2}  & =\tfrac\partial{\partial x_{1}^{2}},\ldots,\kappa_{m}%
^{2}=\tfrac\partial{\partial x_{m}^{2}}\\ \kappa_{0}^{2}  &
=x_{1}^{2}\tfrac\partial{\partial x_{1}^{1}}+\cdots
+x_{m}^{2}\tfrac\partial{\partial x_{m}^{1}}+x_{1}^{1}\tfrac\partial{\partial
x_{1}^{0}}+\cdots+x_{m}^{1}\tfrac\partial{\partial x_{m}^{0}}+\tfrac
\partial{\partial x_{0}^{0}}.
\end{align*}
Loosely speaking, we can write
\begin{align*}
\kappa_{1}^{2}  & =\tfrac\partial{\partial x_{1}^{2}},\ldots,\kappa_{m}%
^{2}=\tfrac\partial{\partial x_{m}^{2}}\\
\kappa_{0}^{2}  & =x_{1}^{2}\kappa_{1}^{1}+\cdots+x_{m}^{2}\kappa_{m}%
^{1}+\kappa_{0}^{1}.
\end{align*}
In order to make this precise we will adopt the following natural notation.
Consider an arbitrary vector field~$f$ given on $J^{n-1}(\mathbb{R}%
,\mathbb{R}^{m})$ by
\[
f=%
{\textstyle\sum\limits_{i=0}^{n-1}}
{\textstyle\sum\limits_{j=1}^{m}}
f_{j}^{i}(\overline{x}^{n-1})\tfrac\partial{\partial x_{j}^{i}}+f_{0}%
^{0}(\overline{x}^{n-1})\tfrac\partial{\partial x_{0}^{0}},
\]
where~$\overline{x}^{n-1}$ denotes the coordinates $x_{0}^{0},x_{1}^{0}%
,\ldots,x_{m}^{0},x_{1}^{1},\ldots,x_{m}^{1},\ldots,x_{1}^{n-1},\ldots
,x_{m}^{n-1}$ of $J^{n-1}(\mathbb{R},\mathbb{R}^{m})$. We can \emph{lift} the
vector field $f$ to a vector field on $J^{n}(\mathbb{R},\mathbb{R}^{m})$, which
we also denote by~$f$, by taking
\[
f=%
{\textstyle\sum\limits_{i=0}^{n-1}}
{\textstyle\sum\limits_{j=1}^{m}}
f_{j}^{i}(\overline{x}^{n-1})\tfrac\partial{\partial x_{j}^{i}}+f_{0}%
^{0}(\overline{x}^{n-1})\tfrac\partial{\partial x_{0}^{0}}+0\cdot
\tfrac\partial{\partial x_{1}^{n}}+\cdots+0\cdot\tfrac\partial{\partial
x_{m}^{n}}.
\]
That is, we lift~$f$ by translating it along the directions $\tfrac
\partial{\partial x_{1}^{n}},\ldots,\tfrac\partial{\partial x_{m}^{n}}$.

\begin{notation}%
[lifts of vector fields] From now on, in any expression of the form
$\kappa_0^n=\sum\nolimits_{i=0}^m\alpha_i(x)\kappa_i^{n-1}$, the vector fields
$\kappa_0^{n-1},\ldots,\kappa_m^{n-1}$ should be considered as their above
defined lifts.
\end{notation}

Let $\kappa^{n-1}=(\kappa_{1}^{n-1},\ldots,\kappa_{m}^{n-1},\kappa_{0}^{n-1})$
denote a family of vector fields defined on
$J^{n-1}(\mathbb{R},\mathbb{R}^{m})$. A \emph{regular prolongation}, with a
parameter~$c^{n}$, of~$\kappa^{n-1}$, denoted by
$\kappa^{n}=R_{c^{n}}(\kappa^{n-1})$, is a family of vector fields
$\kappa^{n}=(\kappa_{1}^{n},\ldots,\kappa_{m}^{n},\kappa_{0}^{n})$ defined on
$J^{n}(\mathbb{R},\mathbb{R}^{m})$ by
\begin{equation}%
\begin{array}
[c]{l}%
\kappa_{1}^{n}=\tfrac\partial{\partial x_{1}^{n}},\ldots,\kappa_{m}^{n}%
=\tfrac\partial{\partial x_{m}^{n}}\\
\kappa_{0}^{n}=(x_{1}^{n}+c_{1}^{n})\kappa_{1}^{n-1}+\cdots+(x_{m}^{n}%
+c_{m}^{n})\kappa_{m}^{n-1}+\kappa_{0}^{n-1},
\end{array}
\label{acta-regular-prolongation}%
\end{equation}
where $c^{n}=(c_{1}^{n},\ldots,c_{m}^{n})$ is a vector of$~m$ real constants. A
\emph{singular prolongation}, with a parameter~$c^{n}$, of~$\kappa^{n-1}$,
denoted by $\kappa^{n}=S_{c_{n}}(\kappa^{n-1})$, is a family of vector fields
$\kappa^{n}=(\kappa_{1}^{n},\ldots,\kappa_{m}^{n},\kappa_{0}^{n})$ defined on
$J^{n}(\mathbb{R},\mathbb{R}^{m})$ by
\begin{equation}%
\begin{array}
[c]{l}%
\kappa_{1}^{n}=\tfrac\partial{\partial x_{1}^{n}},\ldots,\kappa_{m}^{n}%
=\tfrac\partial{\partial x_{m}^{n}}\\
\kappa_{0}^{n}=(x_{1}^{n}+c_{1}^{n})\kappa_{1}^{n-1}+\cdots+(x_{m-1}%
^{n}+c_{m-1}^{n})\kappa_{m-1}^{n-1}+\kappa_{m}^{n-1}+x_{m}^{n}\kappa_{0}%
^{n-1},
\end{array}
\label{acta-singular-prolongation}%
\end{equation}
where $c^{n}=(c_{1}^{n},\ldots,c_{m-1}^{n},0)$ is a vector of$~m$ real
constants, the last one being zero.

A family of vector fields $\kappa^{n}$ on $J^{n}(\mathbb{R},\mathbb{R}^{m})$,
for $n\geq1$, will be called an \emph{extended Kumpera-Ruiz normal form} if
$\kappa^{n}=\sigma_{n}\circ\cdots\circ\sigma_{2}(\kappa^{1})$, where for each
$2\leq i\leq n$ the map~$\sigma_{i}$ equals either~$R_{c^{i}}$ or~$S_{c^{i}}$,
for some vector parameters$~c^{i}$. In other words, a Kumpera-Ruiz normal form
is a family of vector fields obtained by successive prolongations from the
family of vector fields that spans the canonical contact system
on~$J^{1}(\mathbb{R},\mathbb{R}^{m})$.

The above defined prolongations and prolongations-based definition of extended
Kumpera-Ruiz normal forms generalizes for contact systems analogous operations
introduced by the authors~\cite{pasillas-respondek-goursat} for Goursat
structures.

Let $x : M \rightarrow \mathbb{R}^{(n+1)m+1} \cong
J^{n}(\mathbb{R},\mathbb{R}^{m})$ be a local coordinate system on a
manifold~$M$, in a neighborhood of a given point~$p$ in~$M$. We will say that
an extended Kumpera-Ruiz normal form on~$J^{n}(\mathbb{R},\mathbb{R}^{m})$,
defined in $x$-coordinates, is \emph{centered}\ at~$p$ if we have $x(p)=0$. For
example, on $J^{2}(\mathbb{R},\mathbb{R}^{2})$, we have the two following
extended Kumpera-Ruiz normal forms
\begin{align*}
& \left(
\begin{array}
[c]{c}%
\tfrac\partial{\partial x_{1}^{2}}%
\end{array}
,
\begin{array}
[c]{c}%
\tfrac\partial{\partial x_{2}^{2}}%
\end{array}
,
\begin{array}
[c]{c}%
x_{1}^{2}\tfrac\partial{\partial x_{1}^{1}}+x_{2}^{2}\tfrac\partial{\partial
x_{2}^{1}}+x_{1}^{1}\tfrac\partial{\partial x_{1}^{0}}+x_{2}^{1}\tfrac
\partial{\partial x_{2}^{0}}+\tfrac\partial{\partial x_{0}^{0}}%
\end{array}
\right) \\ & \left(
\begin{array}
[c]{c}%
\tfrac\partial{\partial x_{1}^{2}}%
\end{array}
,
\begin{array}
[c]{c}%
\tfrac\partial{\partial x_{2}^{2}}%
\end{array}
,
\begin{array}
[c]{c}%
x_{1}^{2}\tfrac\partial{\partial x_{1}^{1}}+\tfrac\partial{\partial x_{2}^{1}%
}+x_{2}^{2}\left(  x_{1}^{1}\tfrac\partial{\partial x_{1}^{0}}+x_{2}^{1}%
\tfrac\partial{\partial x_{2}^{0}}+\tfrac\partial{\partial x_{0}^{0}}\right)
\end{array}
\right)  ,
\end{align*}
defined by $R_{(0,0)}(\kappa^{1})$ and $S_{(0,0)}(\kappa^{1})$, respectively.
These two normal forms are obviously centered at zero.

The following theorem is the second main contribution of the paper. It asserts
that extended Kumpera-Ruiz normal forms serve as local normal forms for all
singular contact systems for curves, that is for all distributions that satisfy
condition~(i) of Theorem~\ref{thm-pfaff-curves} but fail to fulfill the
regularity condition~(ii) of that theorem.

\begin{theorem}[extended Kumpera-Ruiz normal forms]
\label{thm-ext-kumpera-ruiz}A distribution $\mathcal{D}$ of rank $m+1$ on a
manifold$~M$ of dimension $(n+1)m+1$ is equivalent, in a small enough
neighborhood of any point$~p$ in$~M$, to a distribution spanned by an extended
Kumpera-Ruiz normal form, centered at$~p$ and defined on a suitably chosen
neighborhood of zero, if and only if each element$~\mathcal{D}^{(i)}$ of its
derived flag has constant rank $(i+1)m+1$ and contains an involutive
subdistribution$~\mathcal{L}_i\subset\mathcal{D}^{(i)}$ that has constant
corank one in$~\mathcal{D}^{(i)}$, for $0\leq i\leq n$.
\end{theorem}

\section{Proof of the\ Two Main Theorems}
\label{sec-two-main}

In this Section we provide a proof of Theorem~\ref{thm-ext-kumpera-ruiz}. We
start with several Lemmas --- which will be used in the proof --- that describe
the geometry of incidence between characteristic
distributions~$\mathcal{C}_{i}$ and involutive corank one
subdistributions~$\mathcal{L}_{i}$ of~$\mathcal{D}^{(i)}$.\ Then we prove
Theorem~\ref{thm-ext-kumpera-ruiz}.\ Finally, we conclude
Theorem~\ref{thm-pfaff-curves} as a corollary of
Theorem~\ref{thm-ext-kumpera-ruiz}.

For $i\geq0$, we will denote by $\mathcal{C}_{i}$ the characteristic
distribution of $\mathcal{D}^{(i)}$. It follows directly from the Jacobi
identity that $\mathcal{C}_{i}\subset\mathcal{C}_{i+1}$. Define $d_{i}%
(p)=\dim\mathcal{D}^{(i)}(p)$ and $c_{i}(p)=\dim\mathcal{C}_{i}(p)$. Moreover,
denote $r_{i}(p)=d_{i+1}(p)-d_{i}(p)$.

Though the following result is a direct consequence of Bryant's algebraic
Lemma, we will give its proof as a warm up exercice.\ Indeed, the method used
to prove the inclusion $\mathcal{C}_{0}\subset\mathcal{L}_{0}$ is also used in
the proofs of inclusions $\mathcal{L}_{0}\subset\mathcal{L}_{1}$ and
$\mathcal{L}_{0}\subset\mathcal{C}_{1}$, which will be considered later.

\begin{lemma}[$\mathcal{C}_0 \subset\mathcal{L}_0$]
\label{lem-C0inL0}Let $\mathcal{D}$ be a distribution such that $\mathcal{D}%
^{(0)}$ and $\mathcal{D}^{(1)}$ have constant ranks~$d_0$ and $d_1\geq
d_0+1$%
, respectively. If the distribution $\mathcal{D}$ contains an involutive
subdistribution $\mathcal{L}_0\subset\mathcal{D}^{(0)}$ that has constant
corank one in $\mathcal{D}^{(0)}$ then the ranks of $\mathcal{D}%
^{(0)}$ and $%
\mathcal{D}^{(1)}$ satisfy $r_0\leq d_0-1$. Moreover:
\begin{enumerate}
\item The characteristic distribution $\mathcal{C}_0$ satisfies $\mathcal{C}%
_0\subset\mathcal{L}_0$;
\item The rank of $\mathcal{C}_0$ is constant and equal to $d_0-r_0-1$.
\end{enumerate}
\end{lemma}%

\proof
Assume that $\mathcal{D}$ contains an involutive subdistribution
$\mathcal{L}_{0}\subset\mathcal{D}^{(0)}$ of constant corank one. The relation
$r_{0}\leq d_{0}-1$ is obvious.\ In order to prove by contradiction that
$\mathcal{C}_{0}\subset\mathcal{L}_{0}$, assume that for some point $p$ the
vector space $\mathcal{C}_{0}(p)$ is not contained in $\mathcal{L}_{0}(p)$.
Then, in a small enough neighborhood of$~p$, we can assume that the
distribution $\mathcal{D}^{(0)}$ is a direct sum $\mathcal{D}^{(0)}%
=(h)\oplus\mathcal{L}_{0}$, where $h$ is a vector field that belongs
to~$\mathcal{C}_{0}$ but that does not belong to $\mathcal{L}_{0}$. Since
$[\mathcal{L}_{0},\mathcal{L}_{0}]\subset\mathcal{D}^{(0)}$ we have
$\mathcal{D}^{(1)}=\mathcal{D}^{(0)}+[h,\mathcal{L}_{0}]$. But since
$[h,\mathcal{D}^{(0)}]\subset\mathcal{D}^{(0)}$ we have $\mathcal{D}%
^{(1)}=\mathcal{D}^{(0)}$, which is impossible because $r_{0}\geq1$.

Now, let us compute the rank of $\mathcal{C}_{0}$. Since the corank of
$\mathcal{L}_{0}$ in $\mathcal{D}^{(0)}$ equals $1$, we have a local
decomposition $\mathcal{D}^{(0)}=(f_{d_{0}})\oplus\mathcal{L}_{0}$, where
$f_{d_{0}}$ is an arbitrary vector field that belongs to $\mathcal{D}^{(0)}$
but that does not belong to $\mathcal{L}_{0}$. Since $\mathcal{D}%
^{(1)}=\mathcal{D}^{(0)}+[f_{d_{0}},\mathcal{L}_{0}]$ we can find, locally, a
family of vector fields
\[
(f_{1},\ldots,f_{r_{0}},f_{r_{0}+1},\ldots,f_{d_{0}-1})
\]
that span $\mathcal{L}_{0}$ and satisfies
\[
\mathcal{D}^{(1)}=\mathcal{D}^{(0)}\oplus([f_{d_{0}},f_{1}],\ldots,[f_{d_{0}%
},f_{r_{0}}]).
\]
It follows that, for $r_{0}+1\leq i\leq d_{0}-1$, we can find some smooth
functions $\alpha_{ij}$ such that $[f_{d_{0}},f_{i}]=%
{\textstyle\sum\nolimits_{j=1}^{r_{0}}}
\alpha_{ij}[f_{d_{0}},f_{j}]\operatorname*{mod}\mathcal{D}^{(0)}$. On the one
hand, we have $\dim\mathcal{C}_{0}(p)\geq d_{0}-r_{0}-1$, at each point $p$,
because the vector fields $h_{i}=f_{i}-%
{\textstyle\sum\nolimits_{j=1}^{r_{0}}}
\alpha_{ij}f_{j}$, for $r_{0}+1\leq i\leq d_{0}-1$, satisfy $[h_{i}%
,\mathcal{D}^{(0)}]\subset\mathcal{D}^{(0)}$ and are pointwise linearly
independent. But, on the other hand, we have $\dim\mathcal{C}_{0}(p)\leq
d_{0}-r_{0}-1$, at any point $p$, because otherwise we would have
$\dim\mathcal{D}^{(1)}(p)\leq d_{1}-1$. Hence $\dim\mathcal{C}_{0}%
(p)=d_{0}-r_{0}-1$, for each point $p$ in the underlying manifold.%
\endproof

\begin{lemma}[$\mathcal{L}_0 \subset\mathcal{L}_1$]
\label{lem-L0inL1}Let $\mathcal{D}$ be a distribution such that $\mathcal{D}%
^{(0)}$, $\mathcal{D}^{(1)}$, and $\mathcal{D}^{(2)}$ have constant ranks $%
d_0$, $d_1\geq d_0+2$, and $d_2\geq d_1+2$, respectively. Assume that each
distribution $\mathcal{D}^{(i)}$, for $i=0$ and$~1$, contains an involutive
subdistribution $\mathcal{L}_i\subset\mathcal{D}^{(i)}$ that has constant
corank one in $\mathcal{D}^{(i)}$. Then $\mathcal{L}_0\subset\mathcal{L}_1$.
\end{lemma}%

\proof
Assume that there exists a point $p$ such that the vector space$~\mathcal{L}%
_{0}(p)$ is not contained in $\mathcal{L}_{1}(p)$. We will show that this
assumption leads to $d_{2}\leq d_{1}+1$. To start with, observe that on the one
hand we have $\dim\mathcal{D}^{(0)}(q)\cap\mathcal{L}_{1}(q)\leq d_{0}-1$, for
any point $q$ in a small enough neighborhood of$~p$, because
$\mathcal{L}_{0}(p)$ is not contained in$~\mathcal{L}_{1}(p)$; but that on the
other hand we have
\begin{align*}
\dim\mathcal{D}^{(1)}(q)  & \geq\dim\left(  \mathcal{D}^{(0)}(q)\cup
\mathcal{L}_{1}(q)\right) \\
& =\dim\mathcal{D}^{(0)}(q)+\dim\mathcal{L}_{1}(q)-\dim\left(  \mathcal{D}%
^{(0)}(q)\cap\mathcal{L}_{1}(q)\right)  ,
\end{align*}
for any point $q$, which implies that $\dim\left(  \mathcal{D}^{(0)}%
(q)\cap\mathcal{L}_{1}(q)\right)  \geq d_{0}-1$. Therefore, locally around $p$,
we have $\dim\mathcal{D}^{(0)}(q)\cap\mathcal{L}_{1}(q)=d_{0}-1$. An
analogous argument can be applied to the intersection $\mathcal{L}_{0}%
\cap\mathcal{L}_{1}$ in order to show that $\dim\left(  \mathcal{L}_{0}%
(q)\cap\mathcal{L}_{1}(q)\right)  =d_{0}-2$.

The above relations between the ranks of $\mathcal{D}^{(0)}$, $\mathcal{D}%
^{(0)}\cap\mathcal{L}_{1}$, and $\mathcal{L}_{0}\cap\mathcal{L}_{1}$ imply
that we can find a local basis $(f_{1},\ldots,f_{d_{0}})$ of $\mathcal{D}%
^{(0)}$ such that $\mathcal{L}_{0}=(f_{1},\ldots,f_{d_{0}-1})$ and
$\mathcal{L}_{0}\cap\mathcal{L}_{1}=(f_{1},\ldots,f_{d_{0}-2})$ and
$\mathcal{D}^{(0)}\cap\mathcal{L}_{1}=(f_{1},\ldots,f_{d_{0}-2},f_{d_{0}})$.
Moreover, we can assume that $\mathcal{C}_{0}=(f_{1},\ldots,f_{c_{0}})$, where
$c_{0}<d_{0}-1$. Indeed, by Lemma~\ref{lem-C0inL0}, we have $\mathcal{C}%
_{0}\subset\mathcal{L}_{0}$ and $\mathcal{C}_{0}\subset\mathcal{C}_{1}%
\subset\mathcal{L}_{1}$, which obviously implies $\mathcal{C}_{0}%
\subset\mathcal{L}_{0}\cap\mathcal{L}_{1}$.

Since the vector field $f_{d_{0}}$ does not belong to $\mathcal{L}_{0}$, we
have $\mathcal{D}^{(1)}=\mathcal{D}^{(0)}+[f_{d_{0}},\mathcal{L}_{0}]$.
Therefore, the distribution $\mathcal{D}^{(1)}$ admits the following local
decomposition
\[
\mathcal{D}^{(1)}=\mathcal{D}^{(0)}\oplus([f_{d_{0}},f_{c_{0}+1}%
],\ldots,[f_{d_{0}},f_{d_{0}-1}]).
\]
Denote $g_{i}=[f_{d_{0}},f_{i}]$, for $c_{0}+1\leq i\leq d_{0}-1$. Since the
vector field $f_{d_{0}-1}$ does not belong to $\mathcal{L}_{1}$, we have
$\mathcal{D}^{(2)}=\mathcal{D}^{(1)}+[f_{d_{0}-1},\mathcal{D}^{(1)}]$.
Therefore, all vector fields in $\mathcal{D}^{(2)}$ are linear combinations of
those belonging to $\mathcal{D}^{(1)}$ and of the vector fields $[f_{d_{0}%
-1},g_{i}]$, for $c_{0}+1\leq i\leq d_{0}-1$.

Now, observe that $\mathcal{L}_{0}\cap\mathcal{L}_{1}\subset\mathcal{C}_{1}$.
Indeed, the distribution $\mathcal{D}^{(1)}$ admits a local decomposition
$\mathcal{D}^{(1)}=(f_{d_{0}-1})\oplus\mathcal{L}_{1}$. But $[\mathcal{L}%
_{0}\cap\mathcal{L}_{1},\mathcal{L}_{1}]\subset\mathcal{L}_{1}$ and
$[\mathcal{L}_{0}\cap\mathcal{L}_{1},f_{d_{0}-1}]\subset\mathcal{L}_{0}$. Hence
$[\mathcal{L}_{0}\cap\mathcal{L}_{1},\mathcal{D}^{(1)}]\subset
\mathcal{D}^{(1)}$.

We claim that each vector field $[f_{d_{0}-1},g_{i}]$, for $c_{0}+1\leq i\leq
d_{0}-2 $, belongs to $\mathcal{D}^{(1)}$. Indeed, the Jacobi identity gives
\[
\lbrack f_{d_{0}-1},g_{i}]+[f_{d_{0}},[f_{i},f_{d_{0}-1}]]+[f_{i},[f_{d_{0}%
-1},f_{d_{0}}]]=0.
\]
On the one hand, the vector field $[f_{i},[f_{d_{0}-1},f_{d_{0}}]]$ belongs to
$\mathcal{D}^{(1)}$ because the vector field $f_{i}$ belongs to $\mathcal{L}%
_{0}\cap\mathcal{L}_{1}$, which is contained in $\mathcal{C}_{1}$. But on the
other hand, the vector field $[f_{d_{0}},[f_{i},f_{d_{0}-1}]]\ $also belongs
to $\mathcal{D}^{(1)}$ because $[f_{i},f_{d_{0}-1}]$ belongs to $\mathcal{L}%
_{0}$. It follows that the vector field $[f_{d_{0}-1},g_{i}]$ belongs to
$\mathcal{D}^{(1)} $. Hence $\mathcal{D}^{(2)}=\mathcal{D}^{(1)}+[f_{d_{0}%
-1},g_{d_{0}-1}]$, which obviously implies $d_{2}\leq d_{1}+1$. It follows
that we must have $\mathcal{L}_{0}\subset\mathcal{L}_{1}$.%
\endproof

\begin{lemma}[$\mathcal{L}_0 \subset\mathcal{C}_1$]
\label{lem-L0inC1}Let $\mathcal{D}$ be a distribution such that $\mathcal{D}%
^{(0)}$, $\mathcal{D}^{(1)}$, and $\mathcal{D}^{(2)}$ have constant ranks $%
d_0$, $d_1\geq d_0+2$, and $d_2\geq d_1+2$, respectively. Assume that each
distribution $\mathcal{D}^{(i)}$, for $i=0$ and$~1$, contains an involutive
subdistribution $\mathcal{L}_i\subset\mathcal{D}^{(i)}$ that has constant
corank one in $\mathcal{D}^{(i)}$. Then $\mathcal{L}_0\subset\mathcal{C}_1$.
\end{lemma}%

\proof
Take local generators $(f_{1},\ldots,f_{d_{0}-1},f_{d_{0}})$ of $\mathcal{D}%
^{(0)}$ such that $$\mathcal{L}_{0}=(f_{1},\ldots,f_{d_{0}-1}).$$ Since by
Lemma~\ref{lem-L0inL1} we have $\mathcal{L}_{0}\subset\mathcal{L}_{1}$, each
vector field $f_{i}$, for $1\leq i\leq d_{0}-1$, belongs to $\mathcal{L}_{1}$.
Now observe that, since the corank of $\mathcal{L}_{1}$ in $\mathcal{D}^{(1)}$
equals $1$, we can take $r_{0}$ vector fields $g_{1},\ldots,g_{r_{0}}$ in
$\mathcal{D}^{(1)}$ that are linearly independent $\operatorname*{mod}%
\mathcal{D}^{(0)}$ and such that each vector field $g_{i}$, for $1\leq i\leq
r_{0}-1$ belongs to $\mathcal{L}_{1}$. It follows that there exist two smooth
functions $\alpha$ and $\beta$ such that
\[
\mathcal{L}_{1}=(f_{1},\ldots,f_{d_{0}-1},g_{1},\ldots,g_{r_{0}-1},\alpha
f_{d_{0}}+\beta g_{r_{0}}).
\]
We want to prove that $\mathcal{L}_{0}\subset\mathcal{C}_{1}$, that is that
$[f_{i},\mathcal{D}^{(1)}]\subset\mathcal{D}^{(1)}$, for $1\leq i\leq d_{0}%
-1$. By the definition of $\mathcal{D}^{(1)}$ we have $[f_{i},\mathcal{D}%
^{(0)}]\subset\mathcal{D}^{(1)}$, for $1\leq i\leq d_{0}-1$, and by the
involutivity of $\mathcal{L}_{1}$ we have $[f_{i},g_{j}]\in\mathcal{D}^{(1)}$,
for $1\leq j\leq r_{0}-1$. Therefore, what remains to prove is that
$[\mathcal{L}_{0},g_{r_{0}}]\subset\mathcal{D}^{(1)}$.

We will prove that $[\mathcal{L}_{0},g_{r_{0}}]\subset\mathcal{D}^{(1)}$ by
contradiction. Assume that, for some $1\leq i\leq d_{0}-1$, there exists a
point $p$ such that $[f_{i},g_{r_{0}}](p)\notin\mathcal{D}^{(1)}(p)$. This
implies that $[f_{i},g_{r_{0}}](q)\notin\mathcal{D}^{(1)}(q)$, for each point
$q $ in a small neighborhood $U$ of$~p$. But, since $\mathcal{L}_{1}$ is
involutive, the vector field $[f_{i},\alpha f_{d_{0}}+\beta g_{r_{0}}]$ belongs
to $\mathcal{D}^{(1)}$, which clearly implies that $\beta [f_{i},g_{r_{0}}]$
also belongs to $\mathcal{D}^{(1)}$. Therefore, we must have $\beta(q)=0$, for
each point $q$ in $U$. It follows that, in a small enough neighborhood of $p$,
we have $\mathcal{D}^{(0)}\subset\mathcal{L}_{1}$, which implies that
$\mathcal{D}^{(1)}\subset\mathcal{L}_{1}$ because $\mathcal{L}_{1}$ is
involutive. Since $\mathcal{L}_{1}$ has corank one in
$\mathcal{D}^{(1)}$, this is impossible. Therefore $[\mathcal{L}_{0},g_{r_{0}%
}]\subset\mathcal{D}^{(1)}$, which implies that $\mathcal{L}_{0}%
\subset\mathcal{C}_{1}$.%
\endproof

\begin{lemma}[canonical distribution]
\label{lem-canonical-L0}Let $\mathcal{D}$ be a distribution such that $%
\mathcal{D}^{(0)}$, $\mathcal{D}^{(1)}$, and $\mathcal{D}^{(2)}$ have
constant ranks $d_0=m+1$, $d_1=2m+1$, and $d_2=3m+1$, respectively. If $%
m\geq2$ then assume, additionally, that each distribution $\mathcal{D}^{(i)}
$, for $i=0$ and $1$, contains an involutive subdistribution $\mathcal{L}%
_i\subset\mathcal{D}^{(i)}$ that has constant corank one in $\mathcal{D}%
^{(i)}$. Under these assumptions, we have $\mathcal{L}_0=\mathcal{C}_1$, that
is the distribution $\mathcal{D}^{(0)}$ contains a unique involutive
subdistribution $\mathcal{L}_0$ that has constant corank one in $\mathcal{D}%
^{(0)}$ and satisfies $[\mathcal{L}_0,\mathcal{D}^{(1)}]\subset\mathcal{D}%
^{(1)}$.
\end{lemma}%

\proof
For $m=1$ the result is well known (see e.g. \cite{canadas-ruiz-length-one},
\cite{kumpera-ruiz}, \cite{martin-rouchon-driftless},
and~\cite{montgomery-zhitomirskii}; see also~\cite{pasillas-respondek-goursat}%
). Now, if $m\geq2$ then Item (ii) of Lemma~\ref{lem-C0inL0} (applied to
$\mathcal{C}_{1}$) implies that $\dim\mathcal{C}_{1}(p)=2(2m+1)-(3m+1)-1=m$ and
Lemma~\ref{lem-L0inC1} that $\mathcal{L}_{0}(p)\subset\mathcal{C}_{1}(p)$, for
each point $p$ in the underlying manifold. But $\dim\mathcal{L}_{0}(p)=m$. Thus
$\mathcal{L}_{0}(p)=\mathcal{C}_{1}(p)$, which implies that
$\mathcal{L}_{0}$ is uniquely characterized by $[\mathcal{L}_{0}%
,\mathcal{D}^{(1)}]\subset\mathcal{D}^{(1)}$.%
\endproof

\medskip\

\noindent The following result is a natural generalization of a theorem used by
E. von Weber~\cite{weber-article} in his study of Goursat structures. In fact,
the main idea we will use in our proof of Theorem~\ref{thm-ext-kumpera-ruiz} is
quite close to Weber's original idea. A good introduction to the work of E. von
Weber is Cartan's paper~\cite{cartan-weber}. In~our own
paper~\cite{pasillas-respondek-goursat}, the two main results
of~\cite{weber-article} are given in a more modern language.

\begin{lemma}[extended Weber normal form]
\label{lem-ext-weber}Let $\mathcal{D}$ be a distribution defined on a
manifold$~M$ of dimension $(n+1)m+1$. Assume that$~\mathcal{D}^{(0)}$ and$~%
\mathcal{D}^{(1)}$ have constant ranks $m+1$ and $2m+1$, respectively, and
that $\mathcal{D}^{(0)}$ contains an involutive subdistribution $\mathcal{L}%
_0\subset\mathcal{D}^{(0)}$ that has constant corank one in $\mathcal{D}%
^{(0)}$ and satisfies $\mathcal{L}_0\subset\mathcal{C}_1$. Then, in a small
enough neighborhood of any point$~p$ in$~M$, the distribution$~\mathcal{D}$ is
equivalent to a distribution spanned on $J^n(\mathbb{R},\mathbb{R}^m)$ by a
family of vector fields that has the following form:
\[
\left(
\begin{array}{c}
\tfrac\partial{\partial y_1^n}
\end{array}
,\ldots,
\begin{array}{c}
\tfrac\partial{\partial y_m^n}
\end{array}
,
\begin{array}{c}
y_1^n\zeta_1^{n-1}+\cdots+y_m^n\zeta_m^{n-1}+\zeta_0^{n-1}
\end{array}
\right) ,
\]
where $\mathcal{L}_0=(\tfrac\partial{\partial y_1^n},\ldots,\tfrac
\partial{\partial y_m^n})$ and the vector fields $\zeta_1^{n-1},\ldots
,\zeta_m^{n-1},\zeta_0^{n-1}$ are lifts of vector fields on $J^{n-1}(\mathbb{R%
},\mathbb{R}^m)$, that is,
\[
\zeta_i^{n-1}=\zeta_i^{n-1}(y_0^0,y_1^0,\ldots,y_m^0,\ldots
,y_1^{n-1},\ldots,y_m^{n-1})\text{, for }0\leq i\leq m\text{.}
\]
Moreover, the set of local coordinates $(y_0^0,y_1^0,\ldots,y_m^0,\ldots
,y_1^n,\ldots,y_m^n)$, from $M$ into $J^n(\mathbb{R},\mathbb{R}^m)$ can be
taken to be centered at$~p$.
\end{lemma}%

\proof
It follows directly from Frobenius' theorem, applied to the
distribution~$\mathcal{L}_{0}$, that the distribution$~\mathcal{D}$ is locally
equivalent to a distribution spanned on$~\mathbb{R}^{(n+1)m+1}$ by a family of
vector fields that has the following form:
\[
\left(
\begin{array}
[c]{c}%
\tfrac\partial{\partial z_{1}^{n}}%
\end{array}
,\ldots,
\begin{array}
[c]{c}%
\tfrac\partial{\partial z_{m}^{n}}%
\end{array}
,
\begin{array}
[c]{l}%
{\textstyle\sum\limits_{i=0}^{n-1}}
{\textstyle\sum\limits_{j=1}^{m}}
\alpha_{j}^{i}(z)\tfrac\partial{\partial z_{j}^{i}}+\tfrac\partial{\partial
z_{0}^{0}}%
\end{array}
\right)  ,
\]
where $\mathcal{L}_{0}=(\tfrac\partial{\partial z_{1}^{n}},\ldots
,\tfrac\partial{\partial z_{m}^{n}})$ and the local coordinates $z_{0}%
^{0},z_{1}^{0},\ldots,z_{m}^{0},\ldots,z_{1}^{n},\ldots,z_{m}^{n}$ are centered
at$~p$.

Since $\dim\mathcal{D}^{(1)}(p)=2m+1$ we can assume, after a permutation of the
$z$-coordinates, if necessary, that the real $m\times m$ matrix
\[
T=\left(  \tfrac{\partial\alpha_{i}^{n-1}}{\partial z_{j}^{n}}\right)
,\text{\quad for }1\leq i\leq m\text{ and }1\leq j\leq m\text{,}
\]
has full rank $m$, in a small enough neighborhood of zero. We can assume,
moreover, that $\alpha_{i}^{n-1}(0)=0$, for $1\leq i\leq m$. Otherwise,
replace the coordinate $z_{i}^{n-1}$ by $z_{i}^{n-1}-z_{0}^{0}\alpha_{i}%
^{n-1}(0)$. Now, we can define a new set of centered local coordinates
\[
(y_{0}^{0},y_{1}^{0},\ldots,y_{m}^{0},\ldots,y_{1}^{n},\ldots,y_{m}^{n}%
)=\psi(z_{0}^{0},z_{1}^{0},\ldots,z_{m}^{0},\ldots,z_{1}^{n},\ldots,z_{m}%
^{n})
\]
by taking $y_{i}^{n}=\alpha_{i}^{n-1}(z)$, for $1\leq i\leq m$, and by taking
$y_{j}^{i}=z_{j}^{i}$, as the remaining coordinates. Since the matrix $T$ has
rank $m$, this change of coordinates is indeed a local diffeomorphism. Hence,
the distribution$~\mathcal{D}$ is locally equivalent to a distribution spanned,
on a small enough neighborhood of zero, by a pair of vector fields that has the
following form:
\[
\left(
\begin{array}
[c]{c}%
\tfrac\partial{\partial y_{1}^{n}}%
\end{array}
,\ldots,
\begin{array}
[c]{c}%
\tfrac\partial{\partial y_{m}^{n}}%
\end{array}
,
\begin{array}
[c]{l}%
{\textstyle\sum\limits_{j=1}^{m}}
y_{j}^{n}\tfrac\partial{\partial y_{j}^{n-1}}+%
{\textstyle\sum\limits_{i=0}^{n-2}}
{\textstyle\sum\limits_{j=1}^{m}}
\beta_{j}^{i}(y)\tfrac\partial{\partial y_{j}^{i}}+\tfrac\partial{\partial
y_{0}^{0}}%
\end{array}
\right)  .
\]

Since $\mathcal{L}_{0}\subset\mathcal{C}_{1}$ we have $[\mathcal{L}%
_{0},\mathcal{D}^{(1)}]\subset\mathcal{D}^{(1)}$. But this inclusion clearly
implies that, for $1\leq i\leq n-2$ and $1\leq j\leq m$, we have $\partial
^{2}\beta_{j}^{i}/\partial y_{k}^{n}\partial y_{l}^{n}\equiv0$, for $1\leq
k\leq m$ and $1\leq l\leq m$. It follows that all functions $\beta_{j}^{i}$,
for $1\leq i\leq n-2$ and $1\leq j\leq m$, are affine with respect to the
variables $y_{1}^{n},\ldots,y_{m}^{n}$, that is
\[
\beta_{j}^{i}(y)=%
{\textstyle\sum\limits_{k=1}^{m}}
a_{jk}^{i}(\overline{y}^{n-1})y_{k}^{n}+a_{j0}^{i}(\overline{y}^{n-1}),
\]
where $\overline{y}^{n-1}$ denotes the coordinates $y_{0}^{0},y_{1}^{0}%
,\ldots,y_{m}^{0},\ldots,y_{1}^{n-1},\ldots,y_{m}^{n-1}$. Now, define
\[
\zeta_{k}^{n-1}=\tfrac\partial{\partial y_{j}^{n-1}}+%
{\textstyle\sum\limits_{i=1}^{n-2}}
{\textstyle\sum\limits_{j=1}^{m}}
a_{jk}^{i}(\overline{y}^{n-1})\tfrac\partial{\partial y_{j}^{i}},\text{\quad
for }1\leq k\leq m,
\]
and
\[
\zeta_{0}^{n-1}=%
{\textstyle\sum\limits_{i=1}^{n-2}}
{\textstyle\sum\limits_{j=1}^{m}}
a_{j0}^{i}(\overline{y}^{n-1})\tfrac\partial{\partial y_{j}^{i}}%
+\tfrac\partial{\partial y_{0}^{0}}.
\]
This definition shows that $\mathcal{D}$ is locally equivalent to
\[
\left(
\begin{array}
[c]{c}%
\tfrac\partial{\partial y_{1}^{n}}%
\end{array}
,\ldots,
\begin{array}
[c]{c}%
\tfrac\partial{\partial y_{m}^{n}}%
\end{array}
,
\begin{array}
[c]{c}%
y_{1}^{n}\zeta_{1}^{n-1}+\cdots+y_{m}^{n}\zeta_{m}^{n-1}+\zeta_{0}^{n-1}%
\end{array}
\right)  ,
\]
where the vector fields $\zeta_{1}^{n-1},\ldots,\zeta_{m}^{n-1},\zeta
_{0}^{n-1} $ are lifts of vector fields on $J^{n-1}(\mathbb{R},\mathbb{R}%
^{m})$. It follows directly from our construction that $\mathcal{L}%
_{0}=(\tfrac\partial{\partial y_{1}^{n}},\ldots,\tfrac\partial{\partial
y_{m}^{n}})$.%
\endproof

\bigskip\

\noindent\textbf{Proof of Theorem~\ref{thm-ext-kumpera-ruiz}:} We will proceed
by induction on the integer $n\geq1$. For $n=1$, the Theorem is a direct
consequence of Lemma~\ref{lem-ext-weber}. Thus, assume that the Theorem is true
for $n-1\geq1$ and consider a rank $m+1$ distribution $\mathcal{D}$, defined on
a manifold $M$ of dimension $(n+1)m+1$, such that each element
$\mathcal{D}^{(i)}$ of its derived flag has constant rank $(i+1)m+1$ and
contains an involutive subdistribution $\mathcal{L}_{i}\subset\mathcal{D}%
^{(i)}$ that has constant corank one in $\mathcal{D}^{(i)}$, for $0\leq i\leq
n$. Let $p$ be an arbitrary point in$~M$.

By Lemma~\ref{lem-canonical-L0}, the involutive distribution $\mathcal{L}_{0}%
$, which has corank one in $\mathcal{D}^{(0)}$, satisfies $\mathcal{L}%
_{0}\subset\mathcal{C}_{1}$. We can thus apply Lemma~\ref{lem-ext-weber}, which
states that the distribution $\mathcal{D}$ is equivalent, in a small
enough neighborhood of $p$, to a distribution spanned on $J^{n}(\mathbb{R}%
,\mathbb{R}^{m})$ by a family of vector fields~$(\zeta_{1}^{n},\ldots
,\zeta_{m}^{n},\zeta_{0}^{n})$ that has the following form:
\begin{align*}
\zeta_{1}^{n}  & =\tfrac\partial{\partial y_{1}^{n}},\ldots,\zeta_{m}%
^{n}=\tfrac\partial{\partial y_{m}^{n}}\\
\zeta_{0}^{n}  & =y_{1}^{n}\zeta_{1}^{n-1}+\cdots+y_{m}^{n}\zeta_{m}%
^{n-1}+\zeta_{0}^{n-1},\label{ext-weber-normal-form-bis}%
\end{align*}
where the vector fields $\zeta_{1}^{n-1},\ldots,\zeta_{m}^{n-1},\zeta
_{0}^{n-1} $ are lifts of vector fields on $J^{n-1}(\mathbb{R},\mathbb{R}%
^{m})$. In the rest of the proof we will assume that $\mathcal{D}=(\zeta
_{1}^{n},\ldots,\zeta_{m}^{n},\zeta_{0}^{n})$. Note that the $y$-coordinates
are centered at zero.

The aim of the proof will be to construct a local change of coordinates
\[
(x_{0}^{0},x_{1}^{0},\ldots,x_{m}^{0},\ldots,x_{1}^{n}%
,\ldots,x_{m}^{n})=\phi^{n}(y_{0}^{0},y_{1}^{0},\ldots,y_{m}^{0},\ldots
,y_{1}^{n},\ldots,y_{m}^{n}),
\]
a Kumpera-Ruiz normal form $(\kappa_{1}^{n},\ldots,\kappa_{m}^{n},\kappa
_{0}^{n})$ on $J^{n}(\mathbb{R},\mathbb{R}^{m})$, and a smooth map $\mu
^{n}:J^{n}(\mathbb{R},\mathbb{R}^{m})\rightarrow GL(m+1,\mathbb{R})$, given by
an $(m+1)\times(m+1)$ matrix $(\mu_{ij}^{n}(y))$, such that
\begin{equation}
\phi_{*}^{n}(\zeta_{i}^{n})=%
{\textstyle\sum\limits_{j=1}^{m}}
(\mu_{ij}^{n}\circ\psi^{n})\kappa_{j}^{n},\text{\quad for }1\leq i\leq
m,\label{proof-ext-kumpera-ruiz-01}%
\end{equation}
and
\begin{equation}
\phi_{*}^{n}(\zeta_{0}^{n})=%
{\textstyle\sum\limits_{j=0}^{m}}
(\mu_{0j}^{n}\circ\psi^{n})\kappa_{j}^{n},\label{proof-ext-kumpera-ruiz-02}%
\end{equation}
where $\psi^{n}=(\phi^{n})^{-1}$ denotes the inverse of the local
diffeomorphism $\phi^{n}$. Moreover, we will ask the $x$-coordinates to be
centered at zero, that is $\phi^{n}(0)=0$. Observe that we take $\mu_{i0}%
^{n}=0 $, for $1\leq i\leq m$, in order to transform the canonical distribution
of our distribution $\mathcal{D}$, which is given by
$\mathcal{L}_{0}=(\zeta_{1}^{n},\ldots,\zeta_{m}^{n})$ (see
Lemma~\ref{lem-canonical-L0}), into the canonical distribution of the
Kumpera-Ruiz normal form, which is given by $(\kappa_{1}^{n},\ldots,\kappa
_{m}^{n})$.

Recall that the vector fields $\zeta_{1}^{n-1},\ldots,\zeta_{m}^{n-1}%
,\zeta_{0}^{n-1}$ are lifts of vector fields on $J^{n-1}(\mathbb{R}%
,\mathbb{R}^{m})$. If we take $\mathcal{F}=(\zeta_{1}^{n-1},\ldots,\zeta
_{m}^{n-1},\zeta_{0}^{n-1})$ then we will obtain a decomposition
$\mathcal{D}^{(1)}=\mathcal{L}_{0}\oplus\mathcal{F}$ . Since $\mathcal{D}%
^{(1)}$ contains an involutive subdistribution $\mathcal{L}_{1}$ that has
constant corank one in $\mathcal{D}^{(1)}$, it follows directly from the
relation $\mathcal{L}_{0}\subset\mathcal{L}_{1}$ (see Lemma~\ref{lem-L0inL1})
that $\mathcal{F}$ contains an involutive subdistribution that has constant
corank one in $\mathcal{F}$. In fact, it is easy to prove (using the relations
$\mathcal{C}_{i}\subset\mathcal{C}_{i+1}$ and $\mathcal{L}_{i}=\mathcal{C}%
_{i}$) that if $\mathcal{D}$ satisfies the conditions of
Theorem~\ref{thm-ext-kumpera-ruiz} on $J^{n}(\mathbb{R},\mathbb{R}^{m})$ then
$\mathcal{F}$ satisfies the conditions of this Theorem on $J^{n-1}%
(\mathbb{R},\mathbb{R}^{m})$. Now, recall that we have assumed that the Theorem
is true on $J^{n-1}(\mathbb{R},\mathbb{R}^{m})$. The distribution $\mathcal{F}$
is thus locally equivalent to a distribution spanned by a Kumpera-Ruiz normal
form on $J^{n-1}(\mathbb{R},\mathbb{R}^{m})$, centered at zero. It follows that
there exists a local diffeomorphism
\[
(x_{0}^{0},x_{1}^{0},\ldots,x_{m}^{0},\ldots,x_{1}^{n-1}%
,\ldots,x_{m}^{n-1})=\phi^{n-1}(y_{0}^{0},y_{1}^{0},\ldots,y_{m}^{0}%
,\ldots,y_{1}^{n-1},\ldots,y_{m}^{n-1}),
\]
a Kumpera-Ruiz normal form $(\kappa_{1}^{n-1},\ldots,\kappa_{m}^{n-1}%
,\kappa_{0}^{n-1})$ on $J^{n-1}(\mathbb{R},\mathbb{R}^{m})$, and a smooth map
$\mu^{n-1}:J^{n-1}(\mathbb{R},\mathbb{R}^{m})\rightarrow GL(m+1,\mathbb{R})$
such that
\[
\phi_{*}^{n-1}(\zeta_{i}^{n-1})=%
{\textstyle\sum\limits_{j=0}^{m}}
(\mu_{ij}^{n-1}\circ\psi^{n-1})\kappa_{j}^{n-1}\text{,\quad for }0\leq i\leq
m\text{,}
\]
where $\psi^{n-1}=(\phi^{n-1})^{-1}$ denotes the inverse of the local
diffeomorphism $\phi^{n-1}$. Note that we have $\phi^{n-1}(0)=0$.

The following Lemma can be easily proved by a direct computation.

\begin{lemma}[triangular tangent maps]
\label{lem-triangular-trick}Let $\phi^n=(\phi^{n-1},\phi_1^n,\ldots,\phi
_m^n)^{\mathrm{T}}$ be a diffeomorphism of $J^n(\mathbb{R},\mathbb{R}^m)$ such
that its first $nm+1$ components, which are given by $\phi^{n-1}$, depend on
the first $nm+1$ coordinates only. Moreover, let $f$ be a vector field on
$J^n(\mathbb{R},\mathbb{R}^m)$ of the form $f=\alpha f^{n-1}+f_n$, where
$\alpha$ is a smooth function on $J^n(\mathbb{R},\mathbb{R}^m)$, the vector
field $f^{n-1}$ is the lift of a vector field on
$J^{n-1}(\mathbb{R},\mathbb{R}^m)$, and the only non-zero components of $f_n$
are those that multiply $\tfrac\partial {\partial
y_1^n},\ldots,\tfrac\partial{\partial y_m^n}$. Then, we have
\begin{equation}
\phi_{*}^n(f)=(\alpha\circ\psi^n)\phi _{*}^{n-1}(f^{n-1}) + {\textstyle
\sum\limits_{i=1}^{m}} \left( (\mathrm{L}_f\phi_i^n)\circ \psi^n\right)
\tfrac\partial{\partial x{_i^n}}. \label{ext-triangular-trick}
\end{equation}
Note that the vector field $\phi_{*}^{n-1}(f^{n-1})$ is lifted to $J^n(\mathbb{%
R},\mathbb{R}^m)$ and that the coordinates $x_1^n,\ldots,x_m^n$ are those given
by $\phi_1^n,\ldots,\phi_m^n$, respectively.
\end{lemma}

\noindent\emph{Regular case:} If $\mu_{00}^{n-1}(0)\neq0$ then we can complete
$\phi^{n-1}$ to a zero preserving diffeomorphism $\phi^{n}$ of~$J^{n}%
(\mathbb{R},\mathbb{R}^{m})$ by taking $\phi^{n}=(\phi^{n-1},\phi_{1}%
^{n},\ldots,\phi_{m}^{n})^{\mathrm{T}}$, where
\[
\phi_{j}^{n}(y)=\dfrac{%
{\textstyle\sum\nolimits_{i=1}^{m}}
\mu_{ij}^{n-1}y_{i}^{n}+\mu_{0j}^{n-1}}{%
{\textstyle\sum\nolimits_{i=1}^{m}}
\mu_{i0}^{n-1}y_{i}^{n}+\mu_{00}^{n-1}}-\frac{\mu_{0j}^{n-1}(0)}{\mu
_{00}^{n-1}(0)}\text{,\quad for }1\leq j\leq m\text{.}
\]
It is easy to check, using Lemma~\ref{lem-mobius} below, that $\phi^{n}$ is a
local diffeomorphism (because $\mu^{n-1}$ is invertible). In this case, we
define $c_{i}^{n}=(\mu_{0i}^{n-1}/\mu_{00}^{n-1})(0)$ for $1\leq i\leq m$,
$\mu_{ij}^{n}=\mathrm{L}_{\zeta_{i}^{n-1}}\phi_{j}^{n}$ for $0\leq i\leq m$
and $1\leq j\leq m$, and $\mu_{00}^{n}=%
{\textstyle\sum\nolimits_{i=1}^{m}}
\mu_{i0}^{n-1}y_{i}^{n}+\mu_{00}^{n-1}$. Moreover, the Kumpera-Ruiz normal form
$(\kappa_{1}^{n},\ldots,\kappa_{m}^{n},\kappa_{0}^{n})$ is defined to be the
regular prolongation, with parameter $c^{n}=(c_{1}^{n},\ldots,c_{m}^{n})$, of
$(\kappa_{1}^{n-1},\ldots,\kappa_{m}^{n-1},\kappa_{0}^{n-1})$. Let us check
that, in this case, relation~(\ref{proof-ext-kumpera-ruiz-02}) holds. By
relation (\ref{ext-triangular-trick}) we have:
\begin{eqnarray*}
\phi_{*}^{n}(\zeta_{0}^{n})  & = & {\textstyle\sum\limits_{i=1}^{m}}
(y_{i}^{n}\circ\psi^{n})\phi_{*}^{n-1}(\zeta_{i}^{n-1})+\phi_{*}^{n-1}
(\zeta_{0}^{n-1})+{\textstyle\sum\limits_{i=1}^{m}} \left(
(\mathrm{L}_{\zeta_{0}^{n}}\phi_{i}^{n})\circ\psi^{n}\right)
\tfrac\partial{\partial x_{i}^{n}}\\ & = & {\textstyle\sum\limits_{i=1}^{m}}
(y_{i}^{n}\circ\psi^{n})\left({\textstyle\sum\limits_{j=0}^{m}}
(\mu_{ij}^{n-1}\circ\psi^{n-1})\kappa_{j}^{n-1}\right) +
{\textstyle\sum\limits_{j=0}^{m}}
(\mu_{0j}^{n-1}\circ\psi^{n-1})\kappa_{j}^{n-1}\\ & & \mbox{} +
{\textstyle\sum\limits_{i=1}^{m}} \left(  \mu_{0i}^{n}\circ\psi^{n}\right)
\kappa_{i}^{n}\\ & = & {\textstyle\sum\limits_{j=0}^{m}} \left(  \left(
{\textstyle\sum\limits_{i=1}^{m}} \mu_{ij}^{n-1}y_{i}^{n}+\mu_{0j}^{n-1}\right)
\circ\psi^{n}\right) \kappa_{j}^{n-1}+ {\textstyle\sum\limits_{i=1}^{m}} \left(
\mu_{0i}^{n}\circ\psi^{n}\right)  \kappa_{i}^{n}\\ & = & \left(  \left(
{\textstyle\sum\limits_{i=1}^{m}} \mu_{i0}^{n-1}y_{i}^{n}+\mu_{00}^{n-1}\right)
\circ\psi^{n}\right) \\ & & \mbox{} \times\left(
{\textstyle\sum\limits_{j=1}^{m}} \left( \dfrac{
{\textstyle\sum\nolimits_{i=1}^{m}} \mu_{ij}^{n-1}y_{i}^{n}+\mu_{0j}^{n-1}}{
{\textstyle\sum\nolimits_{i=1}^{m}}
\mu_{i0}^{n-1}y_{i}^{n}+\mu_{00}^{n-1}}\circ\psi^{n}\right)  \kappa_{j}
^{n-1}+\kappa_{0}^{n-1}\right)  + {\textstyle\sum\limits_{i=1}^{m}} \left(
\mu_{0i}^{n}\circ\psi^{n}\right)  \kappa_{i}^{n}\\ & = & \left(
\mu_{00}^{n}\circ\psi^{n}\right)  \left( {\textstyle\sum\limits_{j=1}^{m}}
(x_{j}^{n}+c_{j}^{n})\kappa_{j}^{n-1}+\kappa_{0}^{n-1}\right)  +
{\textstyle\sum\limits_{i=1}^{m}} \left(  \mu_{0i}^{n}\circ\psi^{n}\right)
\kappa_{i}^{n}\\ & = & {\textstyle\sum\limits_{i=0}^{m}} \left(
\mu_{0i}^{n}\circ\psi^{n}\right)  \kappa_{i}^{n}.
\end{eqnarray*}
Moreover, for $1\leq i\leq m$, we have
\[
\phi_{*}^{n}(\zeta_{i}^{n})= {\textstyle\sum\limits_{j=1}^{m}} \left(
(\mathrm{L}_{\zeta_{i}^{n}}\phi_{j}^{n})\circ\psi^{n}\right) \kappa_{j}^{n}=
{\textstyle\sum\limits_{j=1}^{m}} \left(  \mu_{ij}^{n}\circ\psi^{n}\right)
\kappa_{j}^{n}.
\]
It follows that both~(\ref{proof-ext-kumpera-ruiz-01})
and~(\ref{proof-ext-kumpera-ruiz-02}) hold.

\medskip\

\noindent\emph{Singular case:} Suppose now that $\mu_{00}^{n-1}(0)=0$. Since
the matrix $\mu^{n-1}$ is invertible in a small enough neighborhood of zero, we
can assume that there exists an integer $1\leq i\leq m$ such that $\mu
_{0i}^{n-1}(0)\neq0$. After a permutation of the coordinates $y_{1}^{n}%
,\ldots,y_{m}^{n}$, if necessary, we can assume that $\mu_{0m}^{n-1}(0)\neq0$.
Now, like in the regular case, we can complete $\phi^{n-1}$ to a zero
preserving diffeomorphism$~\phi^{n}$ of~$J^{n}(\mathbb{R},\mathbb{R}^{m})$ by
taking $\phi^{n}=(\phi^{n-1},\phi_{1}^{n},\ldots,\phi_{m}^{n})^{\mathrm{T}}$,
where
\[
\phi_{j}^{n}(y)=\dfrac{%
{\textstyle\sum\nolimits_{i=1}^{m}}
\mu_{ij}^{n-1}y_{i}^{n}+\mu_{0j}^{n-1}}{%
{\textstyle\sum\nolimits_{i=1}^{m}}
\mu_{im}^{n-1}y_{i}^{n}+\mu_{0m}^{n-1}}-\frac{\mu_{0j}^{n-1}(0)}{\mu
_{0m}^{n-1}(0)}\text{,\quad for }1\leq j\leq m-1\text{,}
\]
and
\[
\phi_{m}^{n}(y)=\dfrac{%
{\textstyle\sum\nolimits_{i=1}^{m}}
\mu_{i0}^{n-1}y_{i}^{n}+\mu_{00}^{n-1}}{%
{\textstyle\sum\nolimits_{i=1}^{m}}
\mu_{im}^{n-1}y_{i}^{n}+\mu_{0m}^{n-1}}.
\]
In this case, we define $c_{i}^{n}=(\mu_{0i}^{n-1}/\mu_{0m}^{n-1})(0)$ for
$1\leq i\leq m-1$. Observe that we can take $c_{m}^{n}=0$ because $\mu
_{00}^{n-1}(0)=0$. We take $\mu_{ij}^{n}=\mathrm{L}_{\zeta_{i}^{n-1}}\phi
_{j}^{n}$, for $0\leq i\leq m$ and $1\leq j\leq m$, and $\mu_{00}^{n}=%
{\textstyle\sum\nolimits_{i=1}^{m}}
\mu_{im}^{n-1}y_{i}^{n}+\mu_{0m}^{n-1}$. Moreover, the Kumpera-Ruiz normal form
$(\kappa_{1}^{n},\ldots,\kappa_{m}^{n},\kappa_{0}^{n})$ is defined to be
the singular prolongation, with parameter $c^{n}=(c_{1}^{n},\ldots,c_{m-1}%
^{n},0)$, of $(\kappa_{1}^{n-1},\ldots,\kappa_{m}^{n-1},\kappa_{0}^{n-1})$. Let
us check that relation~(\ref{proof-ext-kumpera-ruiz-02}) holds. We have:
\begin{eqnarray*}
\phi_{*}^{n}(\zeta_{0}^{n})  & = & {\textstyle\sum\limits_{i=1}^{m}}
(y_{i}^{n}\circ\psi^{n})\phi_{*}^{n-1}(\zeta_{i}^{n-1})+\phi_{*}^{n-1}
(\zeta_{0}^{n-1})+ {\textstyle\sum\limits_{i=1}^{m}} \left(
(\mathrm{L}_{\zeta_{0}^{n}}\phi_{i}^{n})\circ\psi^{n}\right)
\tfrac\partial{\partial x_{i}^{n}}\\ & = & {\textstyle\sum\limits_{j=0}^{m}}
\left(  \left( {\textstyle\sum\limits_{i=1}^{m}}
\mu_{ij}^{n-1}y_{i}^{n}+\mu_{0j}^{n-1}\right) \circ\psi^{n}\right)
\kappa_{j}^{n-1}+ {\textstyle\sum\limits_{i=1}^{m}} \left(
\mu_{0i}^{n}\circ\psi^{n}\right)  \kappa_{i}^{n}\\  & = & \left(
{\textstyle\sum\limits_{i=1}^{m}} \left(
\mu_{im}^{n-1}y_{i}^{n}+\mu_{0m}^{n-1}\right)  \circ\psi^{n}\right) \\ & &
\mbox{} \times\left( {\textstyle\sum\limits_{j=1}^{m-1}} \left(  \dfrac{
{\textstyle\sum\nolimits_{i=1}^{m}} \mu_{ij}^{n-1}y_{i}^{n}+\mu_{0j}^{n-1}}{
{\textstyle\sum\nolimits_{i=1}^{m}}
\mu_{im}^{n-1}y_{i}^{n}+\mu_{0m}^{n-1}}\circ\psi^{n}\right)
\kappa_{j}^{n-1}+\kappa_{m}^{n-1}\right. \\ & & \mbox{} + \left.  \left(
\dfrac{ {\textstyle\sum\nolimits_{i=1}^{m}}
\mu_{i0}^{n-1}y_{i}^{n}+\mu_{00}^{n-1}}{ {\textstyle\sum\nolimits_{i=1}^{m}}
\mu_{im}^{n-1}y_{i}^{n}+\mu_{0m}^{n-1}}\circ\psi^{n}\right)
\kappa_{0}^{n-1}\right)  + {\textstyle\sum\limits_{i=1}^{m}} \left(
\mu_{0i}^{n}\circ\psi^{n}\right)  \kappa_{i}^{n}\\  & = & \left(
\mu_{00}^{n}\circ\psi^{n}\right)  \left( {\textstyle\sum\limits_{j=1}^{m-1}}
(x_{j}^{n}+c_{j}^{n})\kappa_{j}^{n-1}+\kappa_{m}^{n-1}+x_{m}^{n}\kappa
_{0}^{n-1}\right)  + {\textstyle\sum\limits_{i=1}^{m}} \left(
\mu_{0i}^{n}\circ\psi^{n}\right)  \kappa_{i}^{n}\\ & = &
{\textstyle\sum\limits_{i=0}^{m}} \left(  \mu_{0i}^{n}\circ\psi^{n}\right)
\kappa_{i}^{n}.
\end{eqnarray*}
Moreover, like in the regular case, we have
\[
\phi_{*}^{n}(\zeta_{i}^{n})=%
{\textstyle\sum\limits_{j=1}^{m}}
\left(  (\mathrm{L}_{\zeta_{i}^{n}}\phi_{j}^{n})\circ\psi^{n}\right)
\kappa_{j}^{n}=%
{\textstyle\sum\limits_{j=1}^{m}}
\left(  \mu_{ij}^{n}\circ\psi^{n}\right)  \kappa_{j}^{n},
\]
for $1\leq i\leq m$. It follows that
relations~(\ref{proof-ext-kumpera-ruiz-01})
and~(\ref{proof-ext-kumpera-ruiz-02}) hold in both cases.

We have thus proved that the conditions of Theorem~\ref{thm-ext-kumpera-ruiz}
are sufficient for converting a distribution into extended Kumpera-Ruiz normal
form. It is straightforward to check that these conditions are also necessary.%
\endproof

\begin{lemma}[M\"obius transformations]
\label{lem-mobius}Consider a real $n\times n$ matrix $M$ that has the following
form:
\[
M=\left(
\begin{array}{cc}
A & b \\ c & d
\end{array}
\right) ,
\]
where $c$ is a row vector and $b$ a column vector, both of dimension $n-1$, the
real constant $d$ is non-zero, and $A$ is a real $(n-1)\times(n-1)$ matrix. The
linear fractional transformation$~\varphi$, from$~\mathbb{R}^{n-1}$
into$~\mathbb{R}^{n-1}$, defined in a small enough neighborhood of zero by $%
\varphi(x)=\left( Ax+b\right) /\left( cx+d\right) $ is a local diffeomorphism
if and only if the matrix $M$ is invertible.
\end{lemma}%

\proof
We have $\varphi_{*}(0)=\left(  Ad-bc\right)  /d^{2}$ and thus $\det
\varphi_{*}(0)=(1/d^{2})\det(Ad-bc)$. But $\det M=(1/d^{n-2})\det(Ad-bc)$.
Hence $\det\varphi_{*}(0)\neq0$ if and only if $\det M\neq0$.%
\endproof

\bigskip\

\noindent\textbf{Proof of Theorem~\ref{thm-pfaff-curves}:} Let $\mathcal{D}$ be
a distribution of rank $m+1$, defined on a manifold $M$ of dimension
$(n+1)m+1$, that satisfies the conditions of Theorem~\ref{thm-pfaff-curves}. In
particular, the distribution $\mathcal{D}$ satisfies the conditions of
Theorem~\ref{thm-ext-kumpera-ruiz}. Therefore there exists a Kumpera-Ruiz
normal form $(\kappa_{1}^{n},\ldots,\kappa_{m}^{n},\kappa_{0}^{n})$ on
$J^{n}(\mathbb{R},\mathbb{R}^{m})$, defined on a small enough neighborhood of
zero, that is equivalent to the distribution $\mathcal{D}$ considered on a
small enough neighborhood of any point $p$ in $M$. We will thus assume that
$\mathcal{D}=(\kappa_{1}^{n},\ldots,\kappa_{m}^{n},\kappa_{0}^{n})$.

Now, if we exclude Engel's case ($m=1$ and $n=2$), for which
Theorem~\ref{thm-ext-kumpera-ruiz} is well known to be true, it is
straightforward to check that if the sequence of prolongations that defines our
Kumpera-Ruiz normal form contains a singular prolongation then there exists
some integer $2\leq i\leq n$ such that the Lie flag of $\mathcal{D}$ satisfies
\[
\dim\mathcal{D}_{i}(0)<(i+1)m+1.
\]
It thus follows that $(\kappa_{1}^{n},\ldots,\kappa_{m}^{n},\kappa_{0}^{n})$
has necessarily been obtained by a sequence of regular prolongations from the
canonical contact system on $J^{1}(\mathbb{R},\mathbb{R}^{m})$. That~is:
\begin{equation}
\mathcal{D}=\left(
\begin{array}
[c]{l}%
\tfrac\partial{\partial x_{1}^{n}}%
\end{array}
,\ldots,
\begin{array}
[c]{l}%
\tfrac\partial{\partial x_{m}^{n}}%
\end{array}
,
\begin{array}
[c]{l}%
{\textstyle\sum\limits_{i=0}^{n-1}}
{\textstyle\sum\limits_{j=1}^{m}}
\left(  x_{j}^{i+1}+c_{j}^{i+1}\right)  \tfrac\partial{\partial x_{j}^{i}%
}+\tfrac\partial{\partial x_{0}^{0}}%
\end{array}
\right)  .\label{proof-thm-pfaff-curves}%
\end{equation}
What remains to prove now is that the
distribution~(\ref{proof-thm-pfaff-curves}), which is defined on a small enough
neighborhood of zero, is locally equivalent to the canonical contact system on
$J^{n}(\mathbb{R},\mathbb{R}^{m})$, also considered on a small enough
neighborhood of zero. In other words, we have to normalize all constants
$c_{j}^{i+1}$ by making them equal to zero. To this end, observe that the Lie
algebra
\[
\mathfrak{g}=\operatorname*{span}\nolimits_{\mathbb{R}}\{\kappa_{1}^{n}%
,\ldots,\kappa_{m}^{n},\kappa_{1}^{n-1},\ldots,\kappa_{m}^{n-1},\ldots
,\kappa_{1}^{1},\ldots,\kappa_{m}^{1},\kappa_{0}^{n}\},
\]
of dimension $(n+1)m+1$, generated by the vector fields $(\kappa_{1}%
^{n},\ldots,\kappa_{m}^{n},\kappa_{0}^{n})$
defining~(\ref{proof-thm-pfaff-curves}), has the same structure constants
independently of the values of the parameters $c_{j}^{i+1}$. Indeed, the only
non-zero Lie brackets are those given by the relations $[\kappa_{j}^{i}%
,\kappa_{0}^{n}]=\kappa_{j}^{i-1}$, for $1\leq i\leq n$ and $1\leq j\leq m$. By
Cartan's theorem~\cite{cartan-methode-equivalence} on equivalence of frames,
our distribution is locally equivalent to the canonical contact system at zero
(see e.g.~\cite{olver-equivalence} for a modern account on Cartan's
equivalence method).%
\endproof

\appendix

\section{Appendix}
\label{app-A}

\begin{lemma}[Bryant]
Let $\mathcal{D}$ be a distribution such that $\mathcal{D}^{(0)}$ and $%
\mathcal{D}^{(1)}$ have constant ranks $d_0$ and $d_1$, respectively. Put
$r_0=d_1-d_0$. Assume
that the distribution~$\mathcal{D}$ contains a subdistribution $\mathcal{B}%
\subset\mathcal{D}$ that has constant corank one in $\mathcal{D}$ and satisfies
$[\mathcal{B},\mathcal{B}]\subset\mathcal{D}$. If $r_0\geq3$ then $\mathcal{B}$
is involutive.
\end{lemma}%

\proof
Assume that $\mathcal{D}$ contains a subdistribution $\mathcal{B}%
\subset\mathcal{D}$ that has constant corank one in $\mathcal{D}$ and satisfies
$[\mathcal{B},\mathcal{B}]\subset\mathcal{D}$. By Lemma~\ref{lem-C0inL0}, the
rank of the characteristic distribution $\mathcal{C}_{0}$ of
$\mathcal{D}^{(0)}$ is constant and $\mathcal{C}_{0}\subset\mathcal{B}$.
Therefore, there exists a local basis $(f_{1},\ldots,f_{d_{0}})$ of
$\mathcal{D}$ such that
\[
\mathcal{C}_{0}=(f_{1},\ldots,f_{c_{0}})\text{\quad
and\quad}\mathcal{B}=(f_{1},\ldots,f_{d_{0}-1}),
\]
Since $\mathcal{B}$ satisfies $[\mathcal{B},\mathcal{B}]
\subset\mathcal{D}$ we have $\mathcal{D}^{(1)}=\mathcal{D}^{(0)}+[f_{d_{0}%
},\mathcal{B}]$ or, more precisely,
\begin{equation}
\mathcal{D}^{(1)}=\mathcal{D}^{(0)}\oplus([f_{d_{0}},f_{c_{0}+1}],
\ldots,[f_{d_{0}},f_{d_{0}-1}]).\label{proof-bryant-trick}%
\end{equation}
Note that the relation $d_{1}\geq d_{0}+3$ implies that $\operatorname*{card}%
\{c_{0}+1,\ldots,d_{0}-1\}\geq3$.

We want to prove that $[\mathcal{B},\mathcal{B}]\subset\mathcal{B}$. Let us
first prove that $[f_{i},f_{j}]\in\mathcal{B}$, for $c_{0}+1\leq i\leq d_{0}-1$
and $c_{0}+1\leq j\leq d_{0}-1$. Consider an arbitrary triple $f_{i}$, $f_{j}$,
and $f_{k}$ of vector fields such that the indices $i$, $j$, and $k$ are
pairwise different and contained in $\{c_{0}+1,\ldots,d_{0}-1\}$ (it is
important to stress that such a triple exists because $d_{1}\geq d_{0}+3$). It
follows from
relation~(\ref{proof-bryant-trick}) that the three vector fields $[f_{d_{0}%
},f_{i}]$, $[f_{d_{0}},f_{j}]$, and $[f_{d_{0}},f_{k}]$ are linearly
independent $\operatorname*{mod}\mathcal{D}^{(0)}$. Moreover, since
$[\mathcal{B},\mathcal{B}]\subset\mathcal{D}^{(0)}$, there exist three
smooth functions $a$, $b$, and $c$ such that $[f_{i},f_{j}]=af_{d_{0}%
}\operatorname*{mod}\mathcal{B}$, $[f_{j},f_{k}]=bf_{d_{0}}%
\operatorname*{mod}\mathcal{B}$, and $[f_{k},f_{i}]=cf_{d_{0}%
}\operatorname*{mod}\mathcal{B}$. The Jacobi identity gives:
\[
\lbrack f_{i},[f_{j},f_{k}]]+[f_{j},[f_{k},f_{i}]]+[f_{k},[f_{i},f_{j}]]=0,
\]
which implies that $[f_{i},bf_{d_{0}}]+[f_{j},cf_{d_{0}}]+[f_{k},af_{d_{0}}]$
belongs to $\mathcal{D}^{(0)}$, and thus that $b[f_{i},f_{d_{0}}%
]+c[f_{j},f_{d_{0}}]+a[f_{k},f_{d_{0}}]$ belongs to $\mathcal{D}^{(0)}$. The
latter relation implies that $a$, $b$, and $c$ are identically zero because
$[f_{d_{0}},f_{i}]$, $[f_{d_{0}},f_{j}]$, and $[f_{d_{0}},f_{k}]$ are linearly
independent $\operatorname*{mod}\mathcal{D}^{(0)}$. It follows that we have
$[f_{i},f_{j}]\in\mathcal{B}$, for $c_{0}+1\leq i\leq d_{0}-1$ and $c_{0}+1\leq
j\leq d_{0}-1$.

Since $[\mathcal{C}_{0},\mathcal{C}_{0}]=\mathcal{C}_{0}$, what remain to prove
is that $[\mathcal{C}_{0},\mathcal{B}]=\mathcal{B}$. The proof
follows again from the Jacobi identity, applied to any triple $f_{i}$, $f_{j}%
$, and $f_{k}$ of pairwise linearly independent vector fields such that $f_{i}$
belongs to $\mathcal{C}_{0}$ and both $f_{j}$ and $f_{k}$ belong to
$\mathcal{B}$ but do not belong to $\mathcal{C}_{0}$.
\endproof

\section{Appendix}
\label{app-B}

In this appendix we will provide, following Bryant~\cite{bryant-thesis}, a way
to check the conditions of Corollaries~\ref{cor-corank-one-involutive}
and~\ref{cor-slide}. Indeed, we will show how to verify whether or not the
Engel rank equals~$1$, and how to construct explicitly the characteristic
distribution of~${\cal D}$ and --- when it exists --- the unique corank one
subdistribution ${\cal B} \subset {\cal D}$ satisfying $[{\cal B}, {\cal
B}]\subset {\cal D}$.

Consider a distribution ${\cal D}$ of constant rank~$d_0$, defined on a
manifold of dimension~$N$. Let $\omega_{1},\dots,\omega_{s_{0}}$, where
$s_{0}=N-d_{0}$, be differential 1-forms locally spanning ${\cal D}^{\perp}$,
the annihilator of ${\cal D}$, which we denote by $$ {\cal
D}^{\perp}=(\omega_{1},\dots,\omega_{s_{0}}). $$ We will denote by~$\cal I$ the
Pffafian system generated by $\omega_{1},\dots,\omega_{s_{0}}$.

For any form $\omega\in {\cal D}^{\perp}$, we put $$ {\cal W}(\omega)=\{f \in
{\cal D}\ :\ \ f \lrcorner \, d \omega \in  {\cal D}^{\perp}\}. $$ Clearly, the
characteristic distribution ${\cal C}$ of ${\cal D}$ is given by
\begin{equation}
\label{cond-C} {\cal C}={\textstyle \bigcap\limits_{i=1}^{s_{0}}} {\cal
W}(\omega_{i}).
\end{equation}

Now assume that ${\cal D}^{(1)}$ is of constant rank $d_{1}>d_{0}$, that is
$r_{0}\geq 1$, or, equivalently, that the first derived system ${\cal I}^{(1)}$
is of constant rank smaller than $s_{0}$. By a direct calculation we can check
(see e.g.~\cite{bryant-chern-gardner-goldschmidt-griffiths}) that the Engel
rank of the distribution ${\cal D}$, or of the corresponding Pfaffian system
$\cal I$, equals~$1$ at $p$ if and only if
\begin{equation}
\label{cond-E} (d\omega_{i}\wedge d\omega_{j})(p)=0 \mod \cal{I},
\end{equation}
for any $1 \leq i \leq j \leq s_{0}$.

Now let us choose a family of differential 1-forms
$\omega_{1},\dots,\omega_{r_{0}},\omega_{r_{0}+1},\dots,\omega_{s_{0}}$ such
that $({\cal D}^{(0)})^{\perp}=(\omega_{1},\dots,\omega_{s_{0}})$ and $({\cal
D}^{(1)})^{\perp}=(\omega_{r_{0}+1},\dots,\omega_{s_{0}})$. Independently of
the value of $r_{0}\geq 2$, the unique distribution ${\cal B}$ satisfying
$[{\cal B}, {\cal B}]\subset {\cal D}$ is given, as shown by
Bryant~\cite{bryant-thesis}, by
\begin{equation}
\label{cond-B} {\cal B} = {\textstyle \sum\limits_{i=1}^{r_{0}}} {\cal
W}(\omega_{i}).
\end{equation}
In fact, Bryant has also proved that it is enough to take in the above sum only
two terms corresponding to any $1\leq i < j \leq r_{0}$. In order to verify, in
the case $r_{0}=2$, the conditions of Corollary~\ref{cor-corank-one-involutive}
we have additionally to check the involutivity of this explicitly calculable
distribution.

\end{document}